\documentclass[11pt]{article}
\usepackage{amsfonts,amsmath,amssymb}
\usepackage{mathrsfs}
\usepackage[latin1]{inputenc}
\newtheorem{theo}{Theorem}
\newenvironment{Theo}{\begin{theo}\slshape}{\end{theo}}
\newtheorem{Lemm}{Lemma}[section]

\newtheorem{Corol}{Corollary}
\newtheorem{rema}{Remark}
\newenvironment{Rema}{\begin{rema}\normalfont}{\end{rema}}
\newtheorem{defi}{Definition}

\newcommand{\p}{\partial}

\def\qed{\hfill$\square$\par \bigskip}
\newenvironment{Demo}[1]{{\bf Proof #1.~}}{\qed}
\newcommand{\R}{\mathbb{R}}
\newcommand{\C}{\mathbb{C}}
\newcommand{\M}{\mathcal{M}}

\newcommand{\g}{\mathrm{g}}
\newcommand{\cg}{c\mathrm{g}}
\newcommand{\dv}{\, \mathrm{dv}_\mathrm{g}^{n}}
\newcommand{\dvc}{\, \mathrm{dv}_{c\mathrm{g}}^{n}}
\newcommand{\dvv}{\, \mathrm{dv}_\mathrm{g}^{2n-1}}
\newcommand{\ds}{\, \mathrm{d}\sigma_\mathrm{g}^{n-1}}
\newcommand{\dsc}{\, \mathrm{d}\sigma_{c\mathrm{g}}^{n-1}}
\newcommand{\dss}{\, \mathrm{d}\sigma_\mathrm{g}^{2n-2}}
\newcommand{\dd}{\mathrm{d}}
\newcommand{\I}{\mathcal{I}}

\newcommand{\s}{\mathbb{S}}

\newcommand{\norm}[1]{\left\Vert#1\right\Vert}
\newcommand{\abs}[1]{\left\vert#1\right\vert}
\newcommand{\set}[1]{\left\{#1\right\}}
\newcommand{\para}[1]{\left(#1\right)}
\newcommand{\cro}[1]{\left[#1\right]}

\newcommand{\seq}[1]{\left<#1\right>}

\newcommand{\To}{\longrightarrow}
\usepackage{graphicx}
\newcommand{\dive}{\textrm{div}}

\usepackage{times}

\begin{document}
\title{Stability estimates for the anisotropic wave equation from
the Dirichlet-to-Neumann map}
\author{\small \bf Mourad Bellassoued \\
\small University of Carthage, Department of Mathematics\\
\small Faculty of Sciences of Bizerte\\
\small 7021 Jarzouna Bizerte, Tunisia\\
\small mourad.bellassoued@fsb.rnu.tn
\and
\small \bf David Dos Santos Ferreira\thanks{David DSF was partially supported by ANR grant Equa-disp.} \\
\small Universit\'e Paris 13 \\
\small CNRS, UMR 7539 LAGA \\
\small 99, avenue Jean-Baptiste Cl\'ement \\
\small F-93 430 Villetaneuse, France \\
\small ddsf@math.univ-paris13.fr }
\date{}
\maketitle
\begin{abstract}
In this article we seek stability estimates in the inverse problem of determining the potential or the velocity in a wave equation in an anisotropic
medium from measured Neumann boundary observations. This information is enclosed in the dynamical
Dirichlet-to-Neumann map associated to the wave equation. We prove in dimension $n\geq 2$ that the
knowledge of the Dirichlet-to-Neumann map for the
wave equation uniquely determines
the electric potential and we prove H\"older-type stability in
determining the potential. We prove similar results for the determination of velocities close to~1.
\end{abstract}

\tableofcontents

\section{Introduction}

In this paper, we are interested in the following inverse boundary value problem: on a Riemannian manifold with boundary, determine the potential or the
velocity --- i.e. the conformal factor within a conformal class of metrics ---  in a wave equation
from the vibrations measured at the boundary. Let $(\M,\,\g)$ be a compact
Riemannian manifold with boundary $\partial \M$. All manifolds will be assumed smooth (which means $\mathcal{C}^\infty$) and oriented.
We denote by $\Delta_\g$ the Laplace-Beltrami operator associated to the metric $\g$. In local coordinates,
  $$ \g(x)= \sum_{j,k=1}^n \g_{jk}(x) {\rm d}x_j \otimes {\rm d}x_k,$$
$\Delta_\g$ is given by
\begin{equation}\label{1.1}
\Delta_\g=\frac{1}{\sqrt{\det \g}}\sum_{j,k=1}^n\frac{\p}{\p
x_j}\para{\sqrt{\det \g}\,\g^{jk}\frac{\p}{\p x_k}}.
\end{equation}
Here $(\g^{jk})$ is the inverse of the metric $\g$ and $\det
\g=\det(\g_{jk})$. Let us consider the following initial boundary
value problem for the wave equation with bounded (real valued) electric potential $q\in L^\infty(\M)$
\begin{equation}\label{1.2}
\left\{
\begin{array}{llll}
\para{\partial^2_t-\Delta_\g +q(x)}u=0,  & \textrm{in }\; (0,T)\times \M,\cr
\\
u(0,\cdot )=0,\quad\p_tu(0,\cdot)=0 & \textrm{in }\; \M,\cr
\\
u=f, & \textrm{on} \,\,(0,T)\times\p \M,
\end{array}
\right.
\end{equation}
where $f\in H^1((0,T)\times\p\M)$. Denote by $\nu=\nu(x)$ the outer normal to
$\p \M$ at $x\in\p \M$, normalized so that
$\displaystyle\sum_{j,k=1}^n \g^{jk}\nu_j\nu_k=1$. We may define the
dynamical Dirichlet-to-Neumann map $\Lambda_{\g,\,q}$ by
\begin{equation}\label{1.3}
\Lambda_{\g,\,q}f=\sum_{j,k=1}^n\nu_j\g^{jk}\frac{\p u}{\p x_k}\Big|_{(0,T)\times\p\M}.
\end{equation}
It is clear that one cannot hope to uniquely determine the metric
$\g=(\g_{jk})$ from the knowledge of the Dirichlet-to-Neumann map
$\Lambda_{\g,\,q}$. As was noted in \cite{[SU]}, the Dirichlet-to-Neumann map
is invariant under a gauge transformation of the metric~$\g$. Namely, given a diffeomorphism $\Psi:\M\to \M$ such that
$\Psi|_{\p\M}={\rm Id}$ one has $\Lambda_{\Psi^*\g,\,q}=\Lambda_{\g,\,q}$
where $\Psi^*\g$ denotes the pullback of the metric $\g$ under $\Psi$.
The inverse problem should therefore be
formulated modulo the natural gauge invariance. Nevertheless, when the problem is restricted to a conformal class
of metrics, there is no such gauge invariance and the inverse problem now takes the form: knowing $\Lambda_{\cg,q}$, can one determine the
conformal factor $c$ and the potential $q$?
\medskip

Belishev and Kurylev gave an affirmative answer in \cite{[BK]} to the general problem of finding a smooth metric from the Dirichlet-to-Neumann map.
Their approach is based on the boundary control method introduced by Belishev \cite{[1]} and uses in an essential way an unique continuation
property. Unfortunately it seems unlikely that this method would provide stability estimates even under geometric and topological restrictions.
Their method also solves the problem of recovering $\g$ through boundary spectral data.
The boundary control method  gave rise to several refinements of the results of \cite{[BK]}: one can cite for instance \cite{[KL]}, \cite{[KKL]}
and \cite{[AKKLT]}.
\medskip


In this paper, the inverse problem under consideration is whether the knowledge of the
Dirichlet-to-Neumann map $\Lambda_{\g,\,q}$ on the boundary
 uniquely determines the electric potential $q$ (with a fixed metric $\g$) and whether the knowledge of
the Dirichlet-to-Neumann map $\Lambda_{\g}=\Lambda_{\g,0}$ uniquely determines the conformal factor of the metric $\g$
within a conformal class.
From the physical viewpoint, our inverse problem consists in
determining the properties (e.g. a dispersion term) of an
inhomogeneous medium by probing it with disturbances generated on
the boundary. The data are responses of the medium to these
disturbances which are measured on the boundary, and the
goal is to recover the potential $q(x)$ and the velocity $c(x)$ which describes the property of the
medium. Here we assume that the medium is quiet initially, and $f$
is a disturbance which is used to probe the medium. Roughly
speaking, the data is $\p_\nu u$ measured on the boundary for
different choices of $f$.
\medskip

In the Euclidian case ($\g=e$) Rakesh and Symes \cite{[Rakesh-Symes]}, \cite{[Rakesh1]} used complex geometrical optics
solutions concentrating near lines with any direction $\omega\in
\s^{n-1}$ to prove that $\Lambda_{e,q}$ determines $q(x)$ uniquely in the wave equation.
In \cite{[Rakesh-Symes]}, $\Lambda_{e,q}$ gives equivalent information to
the responses on the whole boundary for all the possible input
disturbances.   Ramm and Sj\"ostrand \cite{[Ramm-Sjostrand]} extended the
results in \cite{[Rakesh-Symes]}
to the case of a potential $q$ depending both on space $x$ and time $t$. Isakov \cite{[Isakov1]}
considered the simultaneous determination of a potential
and a damping coefficient.  A key ingredient in the existing results, is the
construction of complex geometric optics solutions of the wave
equation in the Euclidian case, concentrated along a line, and the relationship between the
hyperbolic Dirichlet-to-Neumann map and the X-ray transform plays a
crucial role.
\medskip

Regarding stability estimates, Sun \cite{[Sun]} established in the Euclidean case stability estimates for potentials from the Dirichlet-to-Neumann
map. In \cite{[SU]} and \cite{[SU2]} Stefanov and Uhlmann considered the inverse problem of determining a Riemannian metric
on a Riemannian manifold with boundary from the hyperbolic Dirichlet-to-Neumann map associated to solutions of the wave equation $(\p_t^2-\Delta_\g)u=0$.
A H\"older type of conditional stability estimate was proven in \cite{[SU]} for metrics close enough to the Euclidean metric in $\mathcal{C}^k$, $k\geq 1$
or for generic simple metrics in \cite{[SU2]}.
\medskip

Uniqueness properties for local Dirichlet-to-Neumann maps associated with the wave equation are rather well understood (e.g.,
Belishev \cite{[1]}, Katchlov, Kurylev and Lassas \cite{[KKL]},
Kurylev and Lassas \cite{[KL]}) but stability for such operators is far from being apprehended.
For instance, one may refer to Isakov and Sun \cite{[IS]} where a local Dirichet-to-Neumann map yields a
stability result in determining a coefficient in a subdomain.
As for results involving a finite number of data in the Dirichlet-to-Neumann map, see
Cheng and Nakamura \cite{[CN]}, Rakesh \cite{[Rakesh1]}.
There are quite a few works on Dirichlet-to-Neumann maps, so our references
are far from being complete: see also  Cardoso and Mendoza \cite{[CM]}, Cheng and Yamamoto \cite{[CY]},
 Eskin \cite{[Eskin1]}-\cite{[Eskin2]}-\cite{[Eskin3]}, Hech and Wang \cite{[Hech-Wang]}, Rachele \cite{[Rachele]}, Uhlmann \cite{[Uhlmann]} as related papers.
\medskip

The main goal of this paper is to study the stability of the inverse problem for the dynamical anisotropic wave equation.
The approach that we develop is a dynamical approach.
Our inverse problem corresponds to a formulation with boundary measurements at infinitely many frequencies.
On the other hand, the main methodology for formulations of inverse problems involving a measurement at a fixed frequency,
is based on $L^2$-weighted inequalites called Carleman estimates. For such applications of Carleman inequalities to inverse problems
we refer for instance to Bellassoued \cite{[Bellassoued]}, Isakov \cite{[I1]}. Most papers treat the determination of spatially varying
functions by a single measurement.  As for observability inequalities by means of
Carleman estimates, see \cite{[Bel-Choul]}, \cite{[Bell-Jel-Yama1]}, \cite{[Bell-Jel-Yama2]}.
\medskip

Our proof is inspired by techniques used by Stefanov and Uhlmann~\cite{[SU2]}, and Dos Santos
Ferreira, Kenig, Salo and Uhlmann \cite{[DKSU]}. In the last reference, an uniqueness theorem for an inverse
problem for an elliptic equation is proved following ideas which in turn go back to
the work of Calder\'on \cite{[Calderon]}. The heuristic underlying idea is that one can (at least formally) translate
techniques used in solving the elliptic equation
     $$ \p_{t}^2+\Delta_{\g} $$
(which is the prototype of equations studied in \cite{[DKSU]}) to the case of the wave equation
    $$ \p_{t}^2-\Delta_{\g} $$
by changing $t$ into $it$. Our problem turns out to be somehow easier because we don't need to construct
complex geometrical solutions, but can rely on classical WKB solutions.
\medskip

\subsection{Weak solutions of the wave equation}
Let $(\M,\g)$ be a (smooth) compact Riemannian
manifold with boundary of dimension $n \geq 2$.
We refer to \cite{[Jost]} for the differential calculus of tensor fields on Riemannian manifolds. If we fix local coordinates $x=\cro{x_1,\ldots,x_n}$ and let
$\cro{\frac{\p}{\p x_1},\dots,\frac{\p}{\p x_n}}$ denote the corresponding tangent vector fields, the inner product and the norm on the tangent
space $T_x\M$ are given by
$$
\g(X,Y)=\seq{X,Y}_\g=\sum_{j,k=1}^n\g_{jk}\alpha_j\beta_k,
$$
$$
\abs{X}_\g=\seq{X,X}_\g^{1/2},\qquad  X=\sum_{i=1}^n\alpha_i\frac{\p}{\p x_i},\quad Y=\sum_{i=1}^n
\beta_i\frac{\p}{\p x_i}.
$$
If $f$ is a $\mathcal{C}^1$ function on $\M$, we define the gradient of $f$ as the vector field $\nabla_\g f$ such that
$$
X(f)=\seq{\nabla_\g f,X}_\g
$$
for all vector fields $X$ on $\M$. In local coordinates, we have
\begin{equation}\label{1.4}
\nabla_\g f=\sum_{i,j=1}^n\g^{ij}\frac{\p f}{\p x_i}\frac{\p}{\p x_j}.
\end{equation}
The metric tensor $\g$ induces the Riemannian volume $\dv=\para{\det \g}^{1/2}\dd x_1\wedge\cdots \wedge \dd x_n$. We denote by $L^2(\M)$ the completion
of $\mathcal{C}^\infty(\M)$ with respect to the usual inner product
$$
\seq{f_1,f_2}=\int_\M f_1(x)\overline{f_2(x)} \dv,\qquad  f_1,f_2\in\mathcal{C}^\infty(\M).
$$
The Sobolev space $H^1(\M)$ is the completion of $\mathcal{C}^\infty(\M)$ with respect to the norm $\norm{\,\cdot\,}_{H^1(\M)}$,
$$
\norm{f}^2_{H^1(\M)}=\norm{f}^2_{L^2(\M)}+\norm{\nabla f}^2_{L^2(\M)}.
$$
The normal derivative is given by
\begin{equation}\label{1.5}
\p_\nu u:=\nabla_\g u\cdot\nu=\sum_{j,k=1}^n\g^{jk}\nu_j\frac{\p u}{\p x_k}
\end{equation}
where $\nu$ is the unit outward vector field to $\p \M$.
Moreover, using covariant derivatives (see \cite{[Hebey]}), it is possible to define coordinate invariant norms in $H^k(\M)$, $k\geq 0$.

Let us consider the following initial boundary value problem
for the wave equation
\begin{equation}\label{1.6}
\left\{
\begin{array}{llll}
\para{\partial_t^2-\Delta_\g+q(x)}v(t,x)=F(t,x)  & \textrm{in }\,\,(0,T)\times\M,\cr
\\
v(0,x)=0,\quad\p_tv(0,x)=0 & \textrm{in }\,\,\M,\cr
\\
v(t,x)=0 & \textrm{on } \,\, (0,T)\times\p\M.
\end{array}
\right.
\end{equation}
The following result is well known (see \cite{[Ikawa]}).
\begin{Lemm}\label{L1.1}
Let $T>0$ and $q\in L^\infty(\M)$, suppose that $F\in \mathscr{H}$, with
$\mathscr{H}=L^1(0,T;L^2(\M))$. The unique solution
$v$ of \eqref{1.6} satisfies
 $$v\in {\cal C}^1(0,T;L^2(\M))\cap {\cal C}(0,T;H^1_0(\M))$$
and  the mapping $F\mapsto\p_\nu v$ is linear and continuous
from $\mathscr{H}$ to $L^2((0,T)\times\p\M)$. Furthermore, there is a
constant $C>0$ such that
\begin{equation}\label{1.7}
\norm{\p_tv(t,\cdot)}_{L^2(\M)}+\norm{\nabla v(t,\cdot)}_{L^2(\M)}\leq
C\norm{F}_{L^1(0,T;L^2(\M))},
\end{equation}
\begin{equation}\label{1.9}
\norm{\p_\nu v}_{L^2((0,T)\times\p\M)}\leq C\norm{F}_{\mathscr{H}}.
\end{equation}
\end{Lemm}
A proof of the following lemma may be found for instance in \cite{[Lions-Magenes]}.
\begin{Lemm}\label{L1.2}
Let $f\in H^1((0,T)\times\p\M)$ be a function such that $f(0,x)=0$ for all $x\in\p\M$.
There exists an unique solution
\begin{equation}\label{1.10}
u\in {\cal C}^1(0,T;L^2(\M))\cap {\cal C}(0,T;H^1(\M))
\end{equation}
to the problem (\ref{1.2}). Furthermore,
the map $f\mapsto \p_\nu u$ is linear and continuous from
$H^1((0,T)\times\p\M)$ into $L^2((0,T)\times \p\M)$.
\end{Lemm}
Therefore the Dirichlet-to-Neumann map  $\Lambda_{\g,q}$ defined by (\ref{1.3}) is continuous. We denote by
$\norm{\Lambda_{\g,q}}$ its norm in $ {\cal L}\para{H^1((0,T)\times\p\M);L^2((0,T)\times \p\M)}$.
Our last remark concerns the fact that when $q$ is real valued, the Dirichlet-to-Neumann map is self-adjoint;
more precisely, we have
\begin{align*}
     \Lambda_{\g,q}^* = \Lambda_{\g,\bar{q}}.
\end{align*}
This simple fact will be proven in section \ref{sec:Preliminaries}. We denote
    $$ \Lambda_{\g} = \Lambda_{\g,0} $$
the Dirichlet-to-Neumann map when there is no potential in the wave equation.

\subsection{Statement of the main results}
In this section we state the main stability results. Let us begin by
introducing an admissible class of manifolds for which we can prove
uniqueness and stability results in our inverse problem. For this we
need the notion of simple manifolds~\cite{[SU2]}.
\begin{defi}
We say that the Riemannian manifold $(\M,\g)$ (or more shortly that the metric $\g$) is simple, if $\p \M$ is
strictly convex with respect to $\g$, and for any $x\in \M$, the exponential map
$\exp_x:\exp_x^{-1}(\M)\To \M$ is a diffeomorphism.
\end{defi}
Note that if $(\M,\g)$ is simple, one can extend $(\M,\g)$ into another simple manifold $\M_{1}$ such that $\M \Subset \M_{1}$.
\medskip

Let us now introduce the admissible set of potentials $q$ and the admissible set of conformal factors $c$. Let $M_0>0$, $k\geq 1$ and
$\varepsilon>0$ be given. Set
\begin{equation}\label{1.11}
\mathscr{Q}(M_0)=\set{q\in H^{1}(\M),\,\,\norm{q}_{H^1(\M)}\leq M_0},
\end{equation}
and
\begin{multline}\label{1.12}
\mathscr{C}(M_0,k,\varepsilon) = \cr
\set{c\in\mathcal{C}^\infty(\M),\,\,c>0\,\,\textrm{in}\,\overline{\M},\,\,
\norm{1-c}_{\mathcal{C}^1(\M)}\leq\varepsilon,\,\,\norm{c}_{\mathcal{C}^k(\M)}\leq M_0}.
\end{multline}
The main results of this paper can be stated as follows.
\begin{Theo}\label{Th1}
Let $(\M,\g)$ be a simple Riemannian compact manifold with boundary of dimension $n \geq 2$, let $T> \mathrm{Diam}_\g(\M)$, there exist constants
$C > 0$ and $\kappa_1\in (0,1)$ such that for any real valued potentials $q_1,\,q_2\in\mathscr{Q}(M_0)$ such that $q_{1}=q_{2}$ on the boundary $\p \M$, we have
\begin{equation}\label{1.13}
\norm{q_1-q_2}_{L^2(\M)}\leq C
\norm{\Lambda_{\g,q_1}-\Lambda_{\g,q_2}}^{\kappa_1}
\end{equation}
where $C$ depends on $\M$, $T$, $M_0$, $n$, and $s$.
\end{Theo}
As a corollary of Theorem \ref{Th1}, we obtain the following uniqueness result.
\begin{Corol}
Let $(\M,\g)$ be a simple Riemannian compact manifold with boundary of dimension $n \geq 2$, let $T>\mathrm{Diam}_\g(\M)$, let $q_1,\,q_2\in\mathscr{Q}(M_0)$
be real valued potentials such that $q_{1}=q_{2}$ on $\p \M$. Then $\Lambda_{\g,q_1}=\Lambda_{\g,q_2}$ implies
$q_1 = q_2$ everywhere in $\M$.
\end{Corol}
\begin{Theo}\label{Th2}
Let $(\M,\g)$ be a simple Riemannian compact manifold with boundary of dimension $n \geq 2$, let $T>\mathrm{Diam}_\g(\M)$, there exist $k\geq 1$,
$\varepsilon>0$, $0<\kappa_2<1$ and $C>0$ such that for any $c\in\mathscr{C}(M_0,k,\varepsilon)$ with
$c=1$ near the boundary $\p\M$, the following estimate holds true
\begin{equation}\label{1.14}
\norm{1-c}_{L^2(\M)}\leq C\norm{\Lambda_{\g}-\Lambda_{\cg}}^{\kappa_2}
\end{equation}
where $C$ depends on $(\M,\g)$, $M_0$, $n$, $\varepsilon$, $k$ and $s$.
\end{Theo}
As a corollary of Theorem \ref{Th2}, we obtain the following uniqueness result.
\begin{Corol}
Let $(\M,\g)$ be a simple Riemannian compact manifold with boundary of dimension $n \geq 2$, let $T>\mathrm{Diam}_\g(\M)$, there exist $k\geq 1$, $\varepsilon>0$, such that for any $c\in\mathscr{C}(M_0,k,\varepsilon)$ with $c=1$ near
the boundary $\p\M$, we have $\Lambda_{\cg}=\Lambda_{\g}$ implies
$c =1$ everywhere in $\M$.
\end{Corol}

\subsection{Spectral inverse problem}
For $q\in \mathscr{Q}(M_0)$ and $q\geq 0$, we denote by $A_q$ the unbounded operator $A_q=-\Delta_\g+q$ with
domain $\mathscr{D}(A_q)=H_0^1(\M )\cap H^2(\M )$.
\medskip

The spectrum of $A_q$ consists of a sequence of eigenvalues, counted according to their multiplicities:
$$
0 \leq \lambda _{1,q}\leq \lambda _{2,q}\leq \ldots \leq \lambda_{k,q} \leq \ldots $$
with $\lim_{k \to \infty} \lambda_{k,q} = \infty$.
The corresponding eigenfunctions are denoted by $(\phi _{k,q})$.
We may assume that this sequence forms an orthonormal basis of
$L^2(\M)$.
\medskip

In the sequel $C$ denotes a generic positive constant depending only on  $\M$ and $M_0$ ($M_0$ is given by (\ref{1.11})).
Since $\phi _{k,q}$ is the solution of the following boundary value problem
\begin{eqnarray*}
\left\{
\begin{array}{ll}
(-\Delta_\g +q)\phi =\lambda _{k,q}\phi  &  \mbox{in}\quad \M \\
\phi=0, & \textrm{on}\quad \p\M,
\end{array}
\right.
\end{eqnarray*}
classical $H^2(\M)$ \textit{a priori} estimates imply
\begin{equation}\label{II.1.1}
\norm{\phi _{k,q}}_{H^\sigma(\M)}\leq C\lambda^{\sigma/2} _{k,q}\norm{
\phi_{k,q}}_{L^2(\M )}=C\lambda _{k,q}^{\sigma/2},\qquad \sigma=0,1,2.
\end{equation}
Therefore
$$
\norm{\partial_\nu \phi _{k,q}}_{H^{1/2}(\p\M)}\leq
 C\lambda _{k,q}.
$$
On the other hand, by Weyl's asymptotics, there exists a positive constant $C\geq 1$ such that
\begin{equation}\label{II.1.2}
C^{-1} k^{2/n}\leq \lambda _{k,q}\leq C k^{2/n}.
\end{equation}
Here $C$ can be chosen uniformly with respect to $q$ provided $0 \le q(x) \le M$ for
$x \in \M$.  Therefore we have
$$
\norm{\partial _\nu \phi _{k,q}}_{H^{1/2}(\p\M )}\leq Ck^{2/n}.
$$
We fix $r$ such that $n/2+1<r \leq n+1$ and it follows
that
$$
\para{k^{-2r/n}\norm{\partial _\nu \phi _{k,q}}_{H^{1/2}(\p\M)}}\in \ell^1.
$$
We recall that $\ell^1$ is the Banach space of real-valued
sequences such that the corresponding series is absolutely
convergent. This space is equipped with its natural norm.
\medskip

Let $\omega=(\omega_k)$ be the sequence given by $\omega_k=k^{-2r/n}$
for each $k\geq 1$. We introduce the following Banach spaces
\begin{multline*}
\ell^1_\omega\para{H^{1/2}(\p\M)}= \cr \quad \set{h=(h_k)_k;\; h_k\in
H^{1/2}\para{\p\M},\; k\geq 1,\; \mbox{and}\; \para{\omega_k\norm{h_k}_{H^{1/2}(\p\M )}}_k\in \ell^1}.
\end{multline*}
and
$$
\ell^1_\omega\para{\C}=\set{y=(y_k)_k;\; y_k\in\C,\; k\geq 1,\; \mbox{and}\; \para{\omega_k \abs{y _k}}_k\in \ell^1}.
$$
The natural norms on those spaces are
$$
\norm{h}_{\ell^1_\omega\para{H^{1/2}(\p\M)}}=\sum_{k\geq
1}\omega_k\norm{h_k}_{H^{1/2}(\p\M)}
$$
and
$$
\norm{y}_{\ell^1_\omega\para{\C}}=\sum_{k\geq
1}\omega_k\abs{y_k}.
$$
We will apply Theorem \ref{Th1} to prove the following result.
\begin{Theo}\label{Th3}
Let $(\M,\g)$ be a simple Riemannian compact manifold with boundary of dimension $n \geq 2$. There exist $C>0$ and $\kappa_3\in(0,1)$  such that the
following estimate holds
\begin{equation}\label{II.1.3}
\norm{q_1-q_2}_{L^2(\M)}\leq C\,\epsilon^{\kappa_3}
\end{equation}
for any non-negative $q_1,q_2\in \mathscr{Q}(M_0)$ which are equal on the boundary $\p\M$, where
$$
\epsilon =\abs{\lambda _{q_1}- \lambda _{q_2}}_{\ell^1_\omega\para{\C}}+\norm{\partial_\nu\phi _{q_1}-\partial _\nu
\phi _{q_2}}_{\ell_\omega^1(H^{1/2}(\p\M))}
$$
is assumed to be small and $\partial _\nu \phi_{q_j}=\para{\partial _\nu\phi _{k,q_j}}_k$, $j=1,2$.
\end{Theo}

Theorem \ref{Th3} is an extension of a result in \cite{[Ch]}
which is itself a variant of a theorem in \cite{[AS]}-\cite{[BCY]}. To the best of our
knowledge, \cite{[AS]} is the first result in the literature concerned with
stability estimates for multidimensional inverse spectral
problems.
\medskip

The outline of the paper is as follows. In section 2 and 3 we collect some  of the formulas needed in the paper. In section 4
we construct special geometrical optics solutions to the wave equation. In section 5 and 6, we establish stability estimates for related
integrals over geodesics crossing $\M$ and prove our main results. In section 7 we prove Theorem 3.

\section{Preliminaries} \label{sec:Preliminaries}
\setcounter{equation}{0}
In this section we collect formulas needed in the rest of this paper.
We denote by $\dive X$ the divergence of a vector field $X\in H^1(T\M)$ on $\M$, i.e. in local coordinates,
\begin{equation}\label{2.1}
\dive X=\frac{1}{\sqrt{\det\g}}\sum_{i=1}^n\p_i\para{\sqrt{\det\g}\,\alpha_i},
\quad X=\sum_{i=1}^n\alpha_i\frac{\p}{\p x_i}.
\end{equation}
If $X\in H^1(T\M)$ we have the divergence formula
\begin{equation}\label{2.2}
\int_\M\dive X \dv=\int_{\p \M}\seq{X,\nu}\ds
\end{equation}
and for $f\in H^1(\M)$ Green's formula reads
\begin{equation}\label{2.3}
\int_\M\dive X \,f\dv=-\int_\M\seq{X,\nabla_\g f}_\g\dv+\int_{\p \M}\seq{X,\nu} f\ds.
\end{equation}
Then if $f\in H^1(\M)$ and $w\in H^2(\M)$, the following identity holds
\begin{equation}\label{2.4}
\int_\M\Delta_\g w f\dv=-\int_\M\seq{\nabla_\g w,\nabla_\g f}_\g\dv+\int_{\p \M}\p_\nu w f \ds.
\end{equation}
Let $f_{1},f_{2} \in H^1((0,T) \times \p \M)$, we denote by $u_{1}$,
respectively by $u_{2}$, the solutions to \eqref{1.2} with potential $q$ and Dirichlet datum $f_{1}$, respectively~$\bar{q}$ and
Dirichlet datum $f_{2}$. By Green's formula, we have
\begin{align*}
     \int_{\p \M} \Lambda_{\g,q}f_{1} \, \overline{f_{2}} \ds &= \int_\M \underbrace{\Delta_\g u_{1} \, \overline{u_{2}}}_{= u_{1} \, \overline{\bar{q}u_{2}}} \dv
     + \int_\M\seq{\nabla_\g u_{1},\overline{\nabla_\g u_{2}}}_\g\dv \\
     &= \int_\M u_{1} \, \overline{\Delta_{\g}u_{2}} \dv + \int_\M\seq{\nabla_\g u_{1},\overline{\nabla_\g u_{2}}}_\g\dv \\
     &= \int_{\p \M} f_{1} \, \overline{\Lambda_{\g,\bar{q}}f_{2}} \ds.
\end{align*}
This shows that
\begin{align*}
     \Lambda_{\g,q}^* = \Lambda_{\g,\bar{q}}.
\end{align*}
In particular, this implies that $\Lambda_{\g,q}$ is selfadjoint when $q$ is real-valued (and therefore $\Lambda_{\g}$). From now on, we will suppose the
potential to be real-valued.

For $x\in \M$ and $\theta\in T_x\M$ we denote by $\gamma_{x,\theta}$ the unique geodesic starting at the point $x$ in the direction $\theta$. We denote
\begin{align*}
S\M&=\set{(x,\theta)\in T\M;\,\abs{\theta}_\g=1},\\
S^*\M&=\set{(x,p)\in T^*\M;\,\abs{p}_\g=1}
\end{align*}
the sphere bundle and co-sphere bundle of $\M$.
The exponential map $\exp_x:T_x\M\To \M$ is given by
\begin{equation}\label{2.5}
\exp_x(v)=\gamma_{x,v}(\abs{v}_\g v)=\gamma_{x,v}(rv),\quad r=\abs{v}_\g.
\end{equation}
A compact Riemannian manifold $(\M,\, \g)$ with boundary is a convex non-trapping
manifold, if it satisfies two conditions:
\begin{enumerate}
 \item[(a)]  the boundary $\p \M$ is strictly
convex, i.e. the second fundamental form of the boundary is positive definite at every boundary point,
 \item[(b)] for every point $x \in \M$ and every vector $\theta\in T_x\M$, $\theta\neq 0$, the maximal
geodesic $\gamma_{x,\theta}(t)$ satisfying the initial conditions
      $$\gamma_{x,\theta}(0) = x \textrm{ and } \dot{\gamma}_{x,\theta}(0) = \theta $$
is defined on a finite segment $[\tau_{-}(x,\theta), \tau_{+}(x,\theta)]$. We recall that a geodesic $\gamma: [a, b] \To M$ is maximal
if it cannot be extended to a segment $[a-\varepsilon_1, b+\varepsilon_2]$, where $\varepsilon_i \geq 0$ and $\varepsilon_1 + \varepsilon_2 > 0$.
\end{enumerate}
The second condition is equivalent to all geodesics having finite length in $\M$.
An important subclass of convex non-trapping manifold are simple manifolds. Recall that a compact Riemannian
manifold $(\M, \g)$ which is simple satisfies the following properties
\begin{enumerate}
   \item[(a)] the boundary is strictly convex,
   \item[(b)] there are no conjugate points on any geodesic.
\end{enumerate}
A simple $n$- dimensional Riemannian manifold is diffeomorphic to a closed ball in $\R^n$, and any pair of points on the manifold can be joined by an unique
minimizing geodesic. \\

In the rest of this article, $C$ will be a generic constant which might change from one line to another, but which only depends on the quantities allowed in the statement of the theorems (namely the quantities involved in the sets $\mathscr{Q}, \mathscr{C}$, the manifold $(\M,g)$, the dimension $n$, the final time $T$ and the H\"older exponents $\kappa_{j}$).

\section{The geodesical ray transform}
\setcounter{equation}{0}
We introduce the submanifolds of inner and outer vectors of $S\M$
\begin{equation}\label{3.1}
\p_{\pm}S\M =\set{(x,\theta)\in S\M,\, x \in \p \M,\, \pm\seq{\theta,\nu(x)}< 0}
\end{equation}
where $\nu$ is the unit outer normal to the boundary. Note that $\p_+ S\M$ and $\p_-S\M$ are
compact manifolds with the same boundary $S(\p \M)$, and $\p S\M = \p_+ S\M\cup \p_- S\M$.
For $(x,\theta)\in\p_+ S\M$,  we denote by $\gamma_{x,\theta} : [0,\tau_+(x,\theta)] \To \M$ the maximal
geodesic satisfying the initial conditions $\gamma_{x,\theta}(0) = x$ and $\dot{\gamma}_{x,\theta}(0) = \theta$.
Let $\mathcal{C}^\infty(\p_+ S\M)$ be the space of smooth functions on the manifold $\p_+S\M$. The ray
transform (also called geodesic X-ray transform) on a convex non-trapping manifold $\M$ is the linear operator
\begin{equation}\label{3.2}
\I:\mathcal{C}^\infty(\M)\To \mathcal{C}^\infty(\p_+S\M)
\end{equation}
defined by
\begin{equation}\label{3.3}
\I f(x,\theta)=\int_0^{\tau_+(x,\theta)}f(\gamma_{x,\theta}(t)) \dd t.
\end{equation}
The right-hand side of (\ref{3.3}) is a
smooth function on $\p_+S\M$ because the integration bound $\tau_+(x,\theta)$ is a smooth function on $\p_+S\M$, see Lemma 4.1.1 of \cite{[Sh]}.
The ray transform on a convex non-trapping manifold $\M$ can be extended to a bounded operator
\begin{equation}\label{3.4}
\I:H^k(\M)\To H^k(\p_+S\M)
\end{equation}
for every integer $k\geq 1$, see Theorem 4.2.1 of \cite{[Sh]}.\\

The Riemannian scalar product on $T_x\M$ induces a volume form on $S_x\M$ denoted by ${\rm d}\omega_x(\theta)$ and given by
$$
{\rm d}\omega_x(\theta)=\sum_{k=1}^n(-1)^k\theta^k {\rm d}\theta^1\wedge\cdots\wedge \widehat{{\rm d}\theta^k}\wedge\cdots\wedge {\rm d}\theta^n.
$$
We introduce the volume form $\dvv$ on the manifold $S\M$
$$
\dvv (x,\theta)=\abs{\dd \omega_x(\theta)\wedge \dv}
$$
where $\dv$ is the Riemannnian volume form on $\M$. By Liouville's theorem, the form $\dvv$ is preserved by the geodesic flow. The
corresponding volume form on the boundary $\p S\M = \set{(x,\theta)\in S\M,\, x\in\p \M}$ is given
by
$$
\dss=\abs{\dd \omega_x(\theta)\wedge \ds}
$$
where $\ds$ is the volume form of $\p \M$.
\medskip

Let $L^2_\mu(\p_+S\M)$ be the space of real valued square integrable functions with respect to the measure $\mu(x,\theta) \dss$ with density
$\mu(x,\theta)=\abs{\seq{\theta,\nu(x)}}$. This Hilbert space is endowed with the scalar product given by
\begin{equation}\label{3.5}
\seq{u,v}_{L^2_\mu(\p_+S\M)}=\int_{\p_+S\M}u(x,\theta) v(x,\theta) \mu(x,\theta)\dss.
\end{equation}
The ray transform $\I$ is a bounded operator from $L^2(\M)$ into
$L^2_\mu(\p_+S\M)$ and its adjoint $\I^*:L^2_\mu(\p_+S\M)\To L^2(\M)$ is given by
\begin{equation}\label{3.6}
\I^*\psi(x)=\int_{S_x\M}\psi^*(x,\theta)\, \dd \omega_x(\theta)
\end{equation}
where $\psi^*$ is the extension of the function $\psi$ from $\p_+S\M$ to $S\M$ constant on
every orbit of the geodesic flow, i.e.
$$
\psi^*(x,\theta)=\psi\big(\gamma_{x,\theta}(\tau_+(x,\theta))\big).
$$
Let $(\M,\g)$ be a simple metric, we assume, as we may, that $(\M,\g)$ extends smoothly into a simple manifold such that $\M_1 \Supset \M$. Then there exist
$C_1>0,C_2>0$ such that
\begin{equation}\label{3.7}
C_1\norm{f}_{L^2(\M)}\leq\norm{\I^*\I(f)}_{H^1(\M_1)}\leq C_2\norm{f}_{L^2(\M)}
\end{equation}
for any $f\in L^2(\M)$, see Theorem 3 in \cite{[SU1]}. If $V$ is an open set of the simple Riemannian manifold $(\M_{1},\g)$, the normal operator $\I^*\I$ is an
elliptic  pseudodifferential operator of order $-1$ on $V$ whose principal symbol is a multiple of $\abs{\xi}_{\g}$ (see \cite{[PU],[SU1]}).
Therefore there exists a constant $C_k>0$ such that for all $f\in H^k(V)$ compactly supported in $V$
\begin{equation}\label{3.8}
\norm{\I^*\I(f)}_{H^{k+1}(\M_{1})}\leq C_k\norm{f}_{H^k(V)}.
\end{equation}

\section{Geometrical optics solutions}
\label{sec:GOSol}
\setcounter{equation}{0}
We will now construct geometrical optics solutions of the  wave equation. We extend the  manifold $(\M,\g)$ into a simple manifold
$\M_2 \Supset \M$ and consider a simple manifold $(\M_1,\g)$ such that $\M_2\Supset \M_1$. The potentials
$q_{1},q_{2}$ may also be extended to $\M_{2}$ and their $H^1(\M_{1})$ norms may be bounded by $M_{0}$. Since $q_{1}$ and $q_{2}$
coincide on the boundary, their extension outside $\M$ can be taken the same so that $q_{1}=q_{2}$ in $\M_{2} \setminus \M_{1}$.
\medskip

Let us assume for a moment that there exist a function $\psi\in{\cal C}^2(\M)$ which satisfies the eikonal equation
\begin{equation}\label{4.1}
\abs{\nabla_\g\psi}^2_\g=\sum_{i,j=1}^n\g^{ij}\frac{\p\psi}{\p x_i}\frac{\p\psi}{\p
x_j}=1,\qquad \forall x\in \M_2
\end{equation}
and a function $a\in H^1(\R,H^2(\M))$ which solves the transport equation
\begin{equation}\label{4.2}
\frac{\p a}{\p t}+\sum_{j,k=1}^n \g^{jk}\frac{\p\psi}{\p x_j}\frac{\p
a}{\p x_k}+\frac{1}{2}(\Delta_\g\psi)a=0,\qquad \forall t\in\R,\, x\in\M
\end{equation}
with initial or final data
\begin{equation}\label{4.3}
a(t,x)=0,\quad \forall x\in \M,\quad \textrm{and}\,\,t\leq0,\,\,\textrm{or}\,\,t\geq T.
\end{equation}
We also introduce the norm $\norm{\cdot}_*$ given by
\begin{equation}\label{4.7}
\norm{a}_*=\norm{a}_{H^1(0,T;H^2(\M))}+\norm{a}_{H^3(0,T;L^2(\M))}.
\end{equation}
\begin{Lemm}\label{L4.1}
Let $q\in L^\infty(\M)$, for any $\lambda>0$, the equation
\begin{eqnarray*}
(\p_t^2-\Delta_\g +q(x))u&=&0,\quad \textrm{in}\quad \M_T:=(0,T)\times\M,\\
u(\kappa,x)&=&\p_tu(\kappa,x)=0,\quad\kappa=0, \textrm{ or } T
\end{eqnarray*}
has a solution of the form
\begin{equation}\label{4.4}
u(t,x)=a(t,x)e^{i\lambda(\psi(x)-t)}+v_\lambda(t,x),
\end{equation}
such that
\begin{equation}\label{4.5}
u\in {\cal C}^1(0,T;L^2(\M))\cap{\cal C}(0,T;H^1(\M)),
\end{equation}
and where $v_\lambda(t,x)$ satisfies
$$
v_\lambda(t,x)=0,\quad\forall (t,x)\in (0,T)\times\p\M ,
$$
$$
v_\lambda(\kappa,x)=0,\quad \p_tv_\lambda(\kappa,x)=0\quad x\in \M,\quad \kappa=0\; \textrm{or}\;\, T
$$
and
\begin{equation}\label{4.6}
\lambda\norm{v_\lambda(t,\cdot)}_{L^2(\M)}+\norm{\p_tv_\lambda(t,\cdot)}_{L^2(\M)}+\norm{\nabla v_\lambda(t,\cdot)}_{L^2(\M)}\leq C\norm{a}_*.
\end{equation}
The constant $C$ depends only on $T$ and $\M$ (that is $C$ does
not depend on $a$ and $\lambda$).
\end{Lemm}
\begin{Demo}{}
We set
\begin{equation}\label{4.9}
k(t,x)=-\para{\partial_t^2-\Delta_\g+q}\para{a(
t,x)e^{i\lambda(\psi-t)}},\quad (t,x)\in(0,T)\times\M.
\end{equation}
To prove our Lemma it would be enough to show that if $v$ solves
\begin{equation}\label{4.8}
\left\{
\begin{array}{ll}
\para{\partial_t^2-\Delta_\g+q}v(t,x)=k(t,x)
& \textrm{in }\,\, (0,T)\times\M,\cr
\\
v(\kappa,x)=0,\quad \p_tv(\kappa,x)=0 & \textrm{in
}\,\,\M,\,\tau=0,\,\textrm{or}\,\,T\cr
\\
v(t,x)=0 & \textrm{on} \,\, (0,T)\times\p\M,
\end{array}
\right.
\end{equation}
then the estimates (\ref{4.6}) holds. We shall prove the estimate
for $\kappa=0$, and the $\kappa=T$ case may be handled in a similar
way. We have
\begin{eqnarray}\label{4.10}
-k(t,x)&=&e^{i\lambda(\psi(x)-
t)}\para{\p_t^2-\Delta_\g+q(x)}\para{a( t,x)}\cr && +2i\lambda
e^{i\lambda(\psi(x)-t)}\para{\p_ta+\sum_{j,k=1}^n
\g^{jk}\frac{\p\psi}{\p x_j}\frac{\p a}{\p
x_k}+\frac{a}{2}\Delta_\g\psi}\cr &&+\lambda^2
a( t,x) e^{i\lambda(\psi(x)-t)}\para{1-\sum_{j,k=1}^n\g^{jk}\frac{\p\psi}{\p x_j}\frac{\p\psi}{\p
x_k}}.
\end{eqnarray}
Taking into account  (\ref{4.1}) and (\ref{4.2}), the right-hand side of (\ref{4.10}) becomes
\begin{equation}\label{4.11}
k(t,x)=-e^{i\lambda(\psi(x)-t)}\para{\p_t^2-\Delta_\g+q}\para{a(t,x)}
\equiv-e^{i\lambda(\psi(x)- t)}k_0(t,x).
\end{equation}
where $k_0\in H^1(0,T;L^2(\M))$ and satisfies
$$
\norm{k_0}_{L^2((0,T)\times\M)}+\norm{\p_tk_0}_{L^2((0,T)\times\M)}\leq C\norm{a}_*.
$$
Since the coefficient $q$ does not depend on $t$, the function
$$
w_\lambda(t,x)=\int_0^tv_\lambda(s,x) \, \dd s
$$
solves the mixed hyperbolic problem (\ref{4.8}) with right-hand side
$$
k_1(t,x)=\int_0^tk(s,x) \dd s=\frac{1}{i\lambda}\int_0^tk_0(s,x)\p_s\para{e^{i\lambda(\psi-s)}} \dd s.
$$
Integrating by parts with respect to $s$, we conclude that
$$
\norm{k_1}_{L^2((0,T)\times\M)}\leq\frac{C}{\lambda}\norm{a}_*.
$$
By Lemma \ref{L1.1}, we find
\begin{equation}\label{4.12}
v_\lambda\in {\cal C}^1(0,T;L^2(\M))\cap{\cal C}(0,T;H^1_0(\M))
\end{equation}
and
\begin{eqnarray}\label{4.13}
\norm{v_\lambda(t,\cdot)}_{L^2(\M)}=\norm{\p_tw_\lambda(t,\cdot)}_{L^2(\M)}\leq
\frac{C}{\lambda}\norm{a}_*.
\end{eqnarray}
Since $\norm{k}_{L^2((0,T)\times\M)}\leq C\norm{a}_*$, using again the energy estimates for the problem (\ref{4.8}), we obtain
\begin{eqnarray}\label{4.14}
\norm{\p_tv_\lambda(t,\cdot)}_{L^2(\M)}+\norm{\nabla v_\lambda (t,\cdot)}_{L^2(\M)}\leq C\norm{a}_*.
\end{eqnarray}
The proof is complete.
\end{Demo}
\begin{Rema}
   In the construction of geometrical optics solutions, it is not necessary to assume that the potential is time independent. In the case where the potential $q$ is also time dependent, one can proceed along the following lines. With the same notations, $w_{\lambda}$ satisfies the equation
    $$(\partial^2_{t}-\Delta_{\g}+q)w_{\lambda} = k_{1} + \int_{0}^t \big(q(t,x)-q(s,x)\big) v_{\lambda}(s,x) \, ds. $$
If one uses Lemma 1.1 on the interval $[0,\tau]$ one gets
\begin{multline*}
    \|\partial_{t}w_{\lambda}(\tau,\cdot)\|_{L^2(\M)} + \|\nabla_{\g}w_{\lambda}(\tau,\cdot)\|_{L^2(\M)} \\ \leq
    \frac{C}{\lambda} \norm{a}_{*} + C \sqrt{T} \|q\|_{L^{\infty}} \int_{0}^{\tau} \|\partial_{t}w(s,\cdot)\|_{L^2(\M)} \, \dd s
\end{multline*}
and Gronwall's inequality allows to conclude
\begin{eqnarray*}
       \|v_{\lambda}(\tau,\cdot)\|_{L^2(\M)} = \|\partial_{t}w_{\lambda}(\tau,\cdot)\|_{L^2(\M)} \leq
       \frac{C}{\lambda} \norm{a}_{*}\Big(1+T\exp\big(C T^{3/2} \|q\|_{L^{\infty}}\big)\Big).
\end{eqnarray*}
\end{Rema}

We now proceed to construct a phase function $\psi$ solution to the eikonal equation (\ref{4.1}) and an amplitude function $a$ solution to the
transport equation (\ref{4.2}).
\medskip

Let $y\in \p \M_1$. Denote points of $\M_1$ by $(r,\theta)$
where $(r,\theta)$ are polar normal coordinates in $\M_1$ with center
$y$. That is $x=\exp_{y}(r\theta)$ where $r>0$ and
$$
\theta\in S_{y}\M_1=\set{\xi\in T_{y}\M_1,\,\,\abs{\xi}_\g=1}.
$$
In these coordinates (which depend on the choice of $y$) the
metric takes the form
$$
\widetilde{\g}(r,\theta)={\rm d}r^2+\g_0(r,\theta)
$$
where $\g_0(r,\theta)$ is a smooth positive definite metric on $S_{y}\M_1$.
For any function $u$ compactly supported in $\M$, we set for $r>0$ and $\theta\in S_y\M_1$
$$
\widetilde{u}(r,\theta)=u(\exp_{y}(r\theta))
$$
where we have extended $u$ by $0$ outside $\M$.
To solve the eikonal equation (\ref{4.1}) it is enough to take
\begin{equation}\label{4.16}
\psi(x)=d_\g(x,y).
\end{equation}
Then by the simplicity assumption, since $y\in \M_2\backslash\overline{\M}$, we have $\psi\in\cal{C}^\infty(\M)$ and
\begin{equation}\label{4.17}
\widetilde{\psi}(r,\theta)=r=d_\g(x,y).
\end{equation}
We now proceed to the transport equation (\ref{4.2}). Recall that
if $f(r)$ is any function of the geodesic distance $r$, then
\begin{equation}\label{4.18}
\Delta_{\widetilde{\g}}f(r)=f''(r)+\frac{\alpha^{-1}}{2}\frac{\p
\alpha}{\p r}f'(r)
\end{equation}
where $\alpha=\alpha(r,\theta)$ denotes the square of the volume element in geodesic polar coordinates.
The transport equation (\ref{4.2}) becomes
\begin{equation}\label{4.19}
\frac{\p \widetilde{a}}{\p t}+\frac{\p \widetilde{\psi}}{\p
r}\frac{\p \widetilde{a}}{\p
r}+\frac{1}{4}\widetilde{a}\alpha^{-1}\frac{\p \alpha}{\p r}\frac{\p
\widetilde{\psi}}{\p r}=0.
\end{equation}
Thus $\widetilde{a}$ satisfy
\begin{equation}\label{4.20}
\frac{\p \widetilde{a}}{\p t}+\frac{\p \widetilde{a}}{\p
r}+\frac{1}{4}\widetilde{a}\alpha^{-1}\frac{\p \alpha}{\p r}=0.
\end{equation}
Let $\phi\in{\cal C}_0^\infty(\R)$ and $b\in H^2(\p_+S\M)$, we choose $\widetilde{a}$ of the form
\begin{equation}\label{4.21}
\widetilde{a}(t,r,\theta)=\alpha^{-1/4}\phi(t-r)b(y,\theta).
\end{equation}
A simple calculation shows that
\begin{equation}\label{4.22}
\frac{\p \widetilde{a}}{\p
t}(t,r,\theta)=\alpha^{-1/4}\phi'(t-r)b(y,\theta).
\end{equation}
and
\begin{equation}\label{4.23}
\frac{\p \widetilde{a}}{\p
r}(t,r,\theta)=-\frac{1}{4}\alpha^{-5/4}\frac{\p\alpha}{\p
r}\phi(t-r)b(y,\theta)-\alpha^{-1/4}\phi'(t-r)b(y,\theta).
\end{equation}
Finally, (\ref{4.23}) and (\ref{4.22}) yield
\begin{equation}\label{4.24}
\frac{\p \widetilde{a}}{\p t}(t,r,\theta)+\frac{\p \widetilde{a}}{\p
r}(t,r,\theta)=-\frac{1}{4}\alpha^{-1}\widetilde{a}(t,r,\theta)\frac{\p\alpha}{\p
r}.
\end{equation}
If we assume that $\mathop{\rm supp}\phi\subset(0,\varepsilon_0)$, with $\varepsilon_0>0$ small enough so that
\begin{equation}
T> {\rm Diam}_\g(\M)+4\varepsilon_0,
\end{equation}
then for any $x=\exp_y(r\theta)\in \M$, it is easy to
see that $\widetilde{a}(t,r,\theta)=0$ if $t\leq 0$ and $t\geq T$.
\begin{Rema}\label{R1}
If $T> {\rm Diam}_\g(\M)+4\varepsilon_0$ and $\cg$ is $\varepsilon$-close to $\g$, then we also have  $T>{\rm  Diam}_{\cg}(\M)+3\varepsilon_0$.
\end{Rema}
\section{Stability estimate for the electric potential}
\setcounter{equation}{0}
In this section, we complete the proof of Theorem \ref{Th1}. We are going to use the geometrical
optics solutions constructed in the previous section; this will provide information on the geodesic ray transform of the difference of electric potentials.
\subsection{Preliminary estimates}
The main purpose of this section is to present a preliminary estimate, which relates the
difference of the potentials to the Dirichlet-to-Neumann map.
As before, we let $q_1,\,q_2\in\mathscr{Q}(M_0)$ be real valued potentials. We set
$$
q =(q_1-q_2) \in H^1_{0}(\M).
$$
Recall that we have extended $q_{1},q_{2}$ as $H^1(\M_{2})$ in such a way that $q=0$ on $\M_{2} \setminus \M$.
\begin{Lemm}\label{L5.1}
There exists $C>0$ such that for any $a_1$, $a_2\in H^1(\R, H^2(\M))$
satisfying the transport equation \eqref{4.2} with initial data \eqref{4.3} the following estimate
holds true:
\begin{multline}\label{5.1}
\abs{\int_{0}^T\!\!\!\!\int_{\M}q(x)a_1(t,x)\overline{a}_2(t,x)\,\dv dt } \leq
\cr
C\para{\lambda^{-1}+\lambda^3\norm{\Lambda_{\g,\,q_1}-\Lambda_{\g,\,q_2}}}\norm{a_1}_*\norm{a_2}_*
\end{multline}
for any sufficiently large $\lambda > 0$.
\end{Lemm}
\begin{Demo}{} First, if $a_2$ satisfies (\ref{4.2}), (\ref{4.3}) and $\lambda$ is sufficiently large, Lemma \ref{L4.1} guarantees
the existence of the  geometrical optics solutions $u_2$
\begin{equation}\label{5.2}
u_2(t,x)=a_2(t,x)e^{i\lambda(\psi(x)-t)}+v_{2,\lambda}(t,x),
\end{equation}
to the equation with the electric potential $q_2$
\begin{eqnarray*}
   \para{\p_t^2-\Delta_{\g}+q_2(x)}u(t,x)=0\quad &&\textrm{in}\,(0,T)\times\M, \\ u(0,\cdot)=0,\quad\p_tu(0,\cdot)=0 \quad &&\textrm{in}\, \M
\end{eqnarray*}
where $v_{2,\lambda}$ satisfies
\begin{equation}\label{5.3}
\begin{array}{lll}
\lambda\norm{v_{2,\lambda}(t,\cdot)}_{L^2(\M)}+\norm{\nabla
v_{2,\lambda}(t,\cdot)}_{L^2(\M)}\leq C\norm{a_2}_*\cr
\\
v_{2,\lambda}(t,x)=0,\quad\forall (t,x)\in\,(0,T)\times\p\M,
\end{array}
\end{equation}
and
$$
u_2\in  {\cal C}^1(0,T;L^2(\M))\cap {\cal C}(0,T;H^1(\M)).
$$
Let us denote by $f_\lambda$ the function
$$
f_\lambda(t,x)=a_2(t,x)e^{i\lambda(\psi(x)-t)},\quad  t\in(0,T),\,\,x\in \p \M.
$$
Let $v$ denote the solution of the following initial boundary value problem
\begin{equation}\label{5.4}
\left\{\begin{array}{lll}
\para{\p_t^2-\Delta_{\g}+q_1} v=0, & (t,x)\in (0,T)\times\M,\cr
\\
v(0,x)=0,\quad\p_tv(0,x)=0, & x\in \M,\cr
\\
v(t,x)=u_2(t,x):=f_{\lambda}(t,x), & (t,x)\in (0,T)\times\p\M.
\end{array}
\right.
\end{equation}
Taking $ w=v-u_2$, one gets
$$
\left\{\begin{array}{lll}
\para{\p_t^2-\Delta_{\g}+q_1(x)}w(t,x)=q(x)u_2(t,x) & (t,x)\in (0,T)\times\M,\cr
\\
w(0,x)=0,\quad\p_tw(0,x)=0, & x\in \M,\cr
\\
w(t,x)=0, & (t,x)\in (0,T)\times\p\M.
\end{array}
\right.
$$
Since $q(x)u_2\in L^1(0,T;L^2(\M))$ by Lemma \ref{L1.1}, we deduce that
$$
w\in {\cal C}^1(0,T;L^2(\M))\cap {\cal C}(0,T;H^1(\M)).
$$
Therefore, we have constructed a particular solution $u_1\in {\cal C}^1(0,T;L^2(\M))\cap {\cal
C}(0,T;H^1(\M))$ to the backward wave equation
\begin{eqnarray*}
\para{\partial_t^2-\Delta_{\g}+q_1(x)}u_1(t,x)=0,  &&\quad (t,x) \in (0,T)\times\M, \\
u_1(T,x)=0,\quad\p_tu_1(T,x)=0  &&\quad x \in \M,
\end{eqnarray*}
having the form
\begin{equation}\label{5.5}
u_1(t,x)=a_1(t,x)e^{i\lambda(\psi(x)-t)}+v_{1,\lambda}(t,x),
\end{equation}
corresponding to the electric potential $q_1$, with
\begin{equation}\label{5.6}
\lambda\norm{v_{1,\lambda}(t,\cdot)}_{L^2(\M)}+\norm{\nabla
v_{1,\lambda}(t,\cdot)}_{L^2(\M)}\leq C\norm{a_1}_*.
\end{equation}
Integrating by parts and using Green's formula (\ref{2.4}), we find
\begin{multline}\label{5.7}
\int_0^T\!\!\!\int_\M\para{\p_t^2-\Delta_\g+q_1(x)}w\overline{u}_1\dv \, \dd t \\
\,=\, \int_0^T\!\!\!\int_\M q(x) u_2\overline{u}_1\dv \, \dd t
\,=\, -\int_0^T\!\!\!\int_{\p \M}\p_\nu w\overline{u}_1\ds \, \dd t.
\end{multline}
Combining (\ref{5.7}) with (\ref{5.4}), we deduce
\begin{align}\label{5.8}
\int_0^T\!\!\!\int_\M q
u_2\overline{u}_1\dv\, \dd t =
-\int_0^T\!\!\!\int_{\p \M}(\Lambda_{\g,\,q_1}-\Lambda_{\g,\,
q_2})(f_{\lambda}) \overline{g}_\lambda \ds\, \dd t
\end{align}
where
$$
g_\lambda(t,x)=a_1(t,x)e^{i\lambda(\psi(x)-t)},\quad (t,x)\in (0,T)\times\p\M.
$$
It follows from (\ref{5.8}), (\ref{5.5}) and (\ref{5.2}) that
\begin{align}\label{5.9}
\int_0^T\!\!\!\int_\M q(x)a_2(t,x)&\overline{a}_1(t,x)\dv\,dt = \cr
& -\int_0^T\!\!\!\int_{\p \M}\para{\Lambda_{\g,\,q_1}-\Lambda_{\g,\,
q_2}}(f_{\lambda})(t,x)\overline{g}_\lambda(t,x) \ds\,\dd t\cr
& -\int_0^T\!\!\!\int_\M q(x)
a_2( t,x)\overline{v}_{1,\lambda}(t,x)e^{i\lambda(\psi-t)}\dv\, \dd t\cr
& -\int_0^T\!\!\!\int_\M q(x)
v_{2,\lambda}(t,x)\overline{a}_1( t,x)e^{-i\lambda(\psi-t)}\dv\, \dd t\cr
& -\int_0^T\!\!\!\int_\M q(x)
v_{2,\lambda}(t,x)\overline{v}_{1,\lambda}(t,x)\dv\, \dd t.
\end{align}
In view of (\ref{5.6}) and (\ref{5.3}), we have
\begin{align}\label{5.10}
\bigg|\int_0^T\!\!\!\int_\M q a_2 \overline{v}_{1,\lambda} e^{i\lambda(\psi-t)}\dv\, \dd t \bigg|
&\leq C\int_0^T\norm{a_2(t,\cdot)}_ {L^2(\M)}\norm{v_{1,\lambda}(t,\cdot)}_{L^2(\M)} \, \dd t \cr
&\leq C\lambda^{-1}\norm{a_2}_*\norm{a_1}_*.
\end{align}
Similarly,  we deduce
\begin{equation*}\label{5.10'}
\abs{\int_0^T\!\!\!\int_\M q(x)
\overline{a}_1( t,x)v_{2,\lambda}(t,x)e^{-i\lambda(\psi-t)}\dv\, \dd t}\leq C\lambda^{-1}\norm{a_1}_*\norm{a_2}_*.
\end{equation*}
Moreover we have
\begin{equation*}\label{5.10"}
\abs{\int_0^T\!\!\!\int_\M q(x)
v_{2,\lambda}(t,x)\overline{v}_{1,\lambda}(t,x)\dv\, \dd t}\leq C\lambda^{-2}\norm{a_1}_*\norm{a_2}_*.
\end{equation*}
On the other hand, by the trace theorem, we find
\begin{align}\label{5.11}
\bigg|\int_0^T\!\!\!\int_{\p \M}&\para{\Lambda_{\g,\,q_1}-\Lambda_{\g,\,
q_2}}(f_{\lambda})\overline{g}_\lambda \ds\, \dd t\big| \cr
& \leq \norm{\Lambda_{\g,\,q_1}-\Lambda_{\g,\, q_2}}
\norm{f_\lambda}_{H^1((0,T)\times\p\M)}\norm{g_\lambda}_{L^2((0,T)\times\p \M)}\cr
& \leq  C\lambda^3\norm{a_1}_*\norm{a_2}_*\norm{\Lambda_{\g,\,q_1}-\Lambda_{\g,\,
q_2}}.
\end{align}
Inequality (\ref{5.1}) follows easily from (\ref{5.9}), (\ref{5.10}), (\ref{5.10'}), (\ref{5.10"}) and (\ref{5.11}).
This completes the proof of the Lemma.
\end{Demo}

\begin{Lemm}\label{L5.2} There exist $C>0$, $\beta\in(0,1)$ such that for any $b\in H^2(\p_+S\M_{1})$, the following estimate
\begin{multline}
\abs{\int_{S_{y}\M_1}\!\int^{\tau_+(y,\theta)}_0  \widetilde{q}(s,\theta) b(y,\theta) \mu(y,\theta) \,\dd s \, \dd \omega_y(\theta)} \\ \leq
C\norm{\Lambda_{\g,\,q_1}-\Lambda_{\g,\,q_2}}^\beta\norm{b(y,\cdot)}_{H^2(S^+_y\M_{1})}
\end{multline}
holds for any $y\in\p\M_1$.
\end{Lemm}
Here we used the notation
    $$ S^+_y\M_{1} = \big\{\theta \in S_{y}\M_{1} : \langle \nu , \theta \rangle < 0 \big\}. $$
We recall that $\mu$ denotes the density $-\langle \theta , \nu(y)\rangle_{\g}$. \\
\begin{Demo}{}
Following (\ref{4.21}), we take two solutions to (\ref{4.2}) and (\ref{4.3}) of the form
\begin{align*}
\widetilde{a}_1(t,r,\theta)&=\alpha^{-1/4}\phi(t-r)b(y,\theta), \\
\widetilde{a}_2(t,r,\theta)&=\alpha^{-1/4}\phi(t-r)\mu(y,\theta).
\end{align*}
Now we change variables in (\ref{5.1}), $x=\exp_{y}(r\theta)$, $r>0$ and
$\theta\in S_{y}\M_1$, we have
\begin{align*}
\int_0^T\!\!\int_\M &q(x)a_1( t,x)a_2(t,x)
\dv\, \dd t\cr
&=\int_0^T\!\!\int_{S_{y}\M_1}\!\int_0^{\tau_+(y,\theta)}\widetilde{q}(r,\theta)\widetilde{a}_1( t,r,\theta)\widetilde{a}_2
(t,r,\theta)
\alpha^{1/2}\, \dd r \, \dd \omega_y(\theta) \, \dd t\cr
&=\int_0^T\!\!\int_{S_{y}\M_1}\!\int_0^{\tau_+(y,\theta)}\widetilde{q}(r,\theta)\phi^2( t-r)b(y,\theta) \mu(y,\theta) \, \dd r \, \dd \omega_y(\theta) \, \dd t.
\end{align*}
By virtue of Lemma \ref{L5.1}, we conclude that
\begin{multline}\label{5.13}
\abs{\int_0^{\infty}\!\!\!\int_{S_{y}\M_1}\!\int_0^{\tau_+(y,\theta)}\widetilde{q}(r,\theta)\phi^2( t-r)b(y,\theta) \mu(y,\theta) \, \dd r \, \dd \omega_y(\theta) \, \dd t} \cr \leq
C\para{\lambda^{-1}+\lambda^3\norm{\Lambda_{\g,\,q_1}-\Lambda_{\g,\,q_2}}}\norm{\phi}^2_{H^3(\R)}\norm{b(y,\cdot)}_{H^2(S^+_y\M_{1})}.
\end{multline}
Since $\phi(t)=0$ for $t\leq 0$ or $t\geq T$, we get
\begin{multline}\label{5.14}
\int_0^{\infty}\!\!\!\int_{S_{y}\M_1}\!\int_0^{\tau_+(y,\theta)}\widetilde{q}(r,\theta)\phi^2( t-r)b(y,\theta) \mu(y,\theta) \, \dd r \, \dd \omega_y(\theta) \, \dd t
\cr = \bigg(\int_{-\infty}^{\infty} \phi^2(t) dt\bigg)  \times \int_{S_{y}\M_1}\!\int_0^{\tau_+(y,\theta)}\widetilde{q}(r,\theta)b(y,\theta) \mu(y,\theta) \, \dd r \, \dd \omega_y(\theta).
\end{multline}
Combining (\ref{5.13}) and (\ref{5.14}), it follows that
\begin{multline*}
 \abs{\int_{S_{y}\M_1}\!\int_0^{\tau_+(y,\theta)}\widetilde{q}(s,\theta)b(y,\theta) \mu(y,\theta) \, \dd s \, \dd\omega_y(\theta)} \\ \leq
C\para{\frac{1}{\lambda}+\lambda^{3}\norm{\Lambda_{\g,\,q_1}-\Lambda_{\g,\,q_2}}}\norm{b(y,\cdot)}_{H^2(S^+_y\M_{1})}.
\end{multline*}
Finally, minimizing in $\lambda$ we obtain
\begin{multline*}
 \abs{\int_{S_{y}\M_1} \! \int_0^{\tau_+(y,\theta)}\widetilde{q}(s,\theta)b(y,\theta) \mu(y,\theta) \, \dd s \, \dd\omega_y(\theta)} \\  \leq
C\norm{\Lambda_{\g,\,q_1}-\Lambda_{\g,\,q_2}}^\beta\norm{b(y,\cdot)}_{H^2(S_y^+\M_{1})}.
\end{multline*}
This completes the proof of the lemma.
\end{Demo}
\subsection{End of the proof of the stability estimate}
Let us now complete the proof of the stability estimate in Theorem \ref{Th1}. Using Lemma \ref{L5.2}, for any $y\in\p \M_1$ and $b\in H^2(\p_+ S\M)$
we have
\begin{multline*}
\abs{\int_{S_{y}\M_1}\I(q)(y,\theta)b(y,\theta)\mu(y,\theta) \, \dd \omega_y(\theta)} \\ \leq C\norm{\Lambda_{\g,\,q_1}-\Lambda_{\g,\,q_2}}^\beta\norm{b(y,\cdot)}_{H^2(S^+_y\M_{1})}.
\end{multline*}
Integrating with respect to $y\in \p \M_1$ we obtain
\begin{multline}\label{5.18}
\abs{\int_{\p_+S\M_1}\I(q)(y,\theta)b(y,\theta)\seq{\theta,\nu(y)}\dss (y,\theta)}\cr  \leq C\norm{\Lambda_{\g,\,q_1}-\Lambda_{\g,\,q_2}}^\beta
\norm{b}_{H^2(\p_+S\M_{1})}.
\end{multline}
Now we choose
$$
b(y,\theta)=\I\para{\I^*\I(q)}(y,\theta).
$$
Taking into account (\ref{3.8}) and (\ref{3.4}), we obtain
$$
\norm{\I^*\I(q)}^2_{L^2(\M_1)}\leq C\norm{\Lambda_{\g,\,q_1}-\Lambda_{\g,\,q_2}}^\beta\norm{q}_{H^1(\M)}.
$$
By interpolation, it follows that
\begin{align}\label{5.19}
\norm{\I^*\I(q)}^2_{H^1(\M_1)}&\leq C \norm{\I^*\I(q)}_{L^2(\M_1)}\norm{\I^*\I(q)}_{H^2(\M_1)}\cr
&\leq C \norm{\I^*\I(q)}_{L^2(\M_1)}\norm{q}_{H^1(\M)}\cr
&\leq C \norm{\I^*\I(q)}_{L^2(\M_1)}\cr
&\leq C\norm{\Lambda_{\g,\,q_1}-\Lambda_{\g,\,q_2}}^{\beta/2}.
\end{align}
Using (\ref{3.7}), we deduce that
$$
\norm{q}_{L^2(\M)}^2\leq C\norm{\Lambda_{\g,\,q_1}-\Lambda_{\g,\,q_2}}^{\beta/2}.
$$
This completes the proof of Theorem \ref{Th1}.
\begin{Rema}
   In the proof of Theorem \ref{Th1}, we have used the time independence of the potential at two stages:
   \begin{enumerate}
    \item In the construction of the remainder term $v_{\lambda}(t,x)$ in the proof of Lemma \ref{L4.1}. But as was noted in Remark \ref{R1}, this
    restriction may be bypassed.
    \item In equation (\ref{5.14}) to get rid of the $\phi^2(t-r)$ term and obtain the ray transform. The adaptation to the time dependent case does
    not seem to be straightforward.
   \end{enumerate}
   We leave the case of a time dependent potential as an open problem.
\end{Rema}

\section{Stability estimate for the conformal factor}
\setcounter{equation}{0}
We shall use the following notations. Let $c\in\mathscr{C}(M_0,k,\varepsilon)$, we denote
\begin{eqnarray}\label{6.1}
&&\varrho_0(x)=1-c(x),\quad\varrho_1(x)=c^{n/2}(x)-1,\quad\varrho_2(x)=c^{n/2-1}(x)-1,\cr
\cr
&&\varrho(x)=\varrho_1(x)-\varrho_2(x)=c^{n/2-1}(x)\para{c(x)-1}.
\end{eqnarray}
Then the following holds
\begin{align}\label{6.2}
\norm{\varrho_j}_{\mathcal{C}^1(\M)}&\leq C\norm{\varrho_0}_{\mathcal{C}^1(\M)},\quad j=1,2.
\\
\label{6.3}
C^{-1}\norm{\varrho_0}_{L^2(\M)}&\leq\norm{\varrho}_{L^2(\M)}\leq C\norm{\varrho_0}_{L^2(\M)}.
\end{align}
The first step in our analysis is the following result.
\begin{Lemm}\label{L6.1}
Let $c\in \mathcal{C}^\infty(\M)$ be such that $c=1$ near the boundary $\p\M$. Let $u_1$, $u_2$ solve the following problem in $(0,T)\times\M$ with
some $T>0$
\begin{eqnarray}\label{6.4}
\left\{
\begin{array}{llll}
(\partial^2_t-\Delta_\g)u_1=0,  & \textrm{in }\; (0,T)\times \M,\cr
\\
u_1(0,\cdot )=\p_tu(0,\cdot)=0 & \textrm{in }\; \M,\cr
\\
u_1=f_1, & \textrm{on} \,\,(0,T)\times\p \M,
\end{array}
\right.
\\ \label{6.4.bis}
\left\{
\begin{array}{llll}
(\partial_t^2-\Delta_{\cg})u_2=0,  & \textrm{in }\; (0,T)\times \M,\cr
\\
u_2(T,\cdot )=\p_tu(T,\cdot)=0 & \textrm{in }\; \M,\cr
\\
u_2=f_2, & \textrm{on} \,\,(0,T)\times\p \M,
\end{array}
\right.
\end{eqnarray}
where $f_k\in H^1\para{(0,T)\times\p\M}$, $k=1,2$. Then the following identity
\begin{multline}\label{6.5} \int_{0}^T\!\!\!\int_{\p\M}\para{\Lambda_\g-\Lambda_{\cg}}f_1\,\overline{f}_2\,\ds \, \dd t
=\int_0^T\!\!\!\int_\M\varrho_1(x)\p_tu_1\p_t\overline{u}_2\,
\dv \, \dd t \\
-\int_0^T\!\!\!\int_\M \varrho_2(x)\seq{\nabla_\g u_1(t,x),\nabla_\g \overline{u}_2(t,x)}_\g\dv \dd t
\end{multline}
holds true for any $f_j\in H^1\para{(0,T)\times\p\M}$, $j=1,2$.
\end{Lemm}
\begin{Demo}{}
We multiply both hand sides of the first equation (\ref{6.4}) by $\overline{u}_2$ ; integrating by parts in time and using Green's formula (\ref{2.3})-(\ref{2.4}) in $(\M,\g)$, we obtain
\begin{align*}
0&=\int_0^T\!\!\!\int_\M\para{\p_t^2u_1-\Delta_\g u_1}\overline{u}_2\dv \dd t \cr
&=-\int_0^T\!\!\!\int_\M \p_t u_1\p_t \overline{u}_2\dvc \, \dd t+\int_0^T\!\!\!\int_\M \varrho_1(x) \p_t u_1\p_t\overline{u}_2 \dv \, \dd t\cr
&\quad+\int_0^T\!\!\!\int_\M\sum_{j,k=1}^n (c\g)^{jk}\para{\frac{\p u_1}{\p x_j}\frac{\p \overline{u}_2}{\p x_k}}\dvc \, \dd t \cr
&\quad -\int_0^T\!\!\!\int_M\varrho_2(x)\para{\sum_{j,k=1}^n
\g^{jk}\frac{\p u_1}{\p x_j}\frac{\p \overline{u}_2}{\p x_k}}\dv \, \dd t
-\int_0^T\!\!\!\int_{\p\M}\p_\nu u_1\overline{f}_2\ds \, \dd t.
\end{align*}
Another integration by parts in time and an application of Green's formula
in $(\M,\cg)$ yield
\begin{align}\label{6.6}
0&=\int_0^T\!\!\!\int_\M u_1\para{\p_t^2\overline{u}_2-\Delta_{\cg}\overline{u}_2}\dvc \, \dd t \cr & \quad +\int_0^T\!\!\!\int_\M \varrho_1(x) \p_tu_1\p_t\overline{u}_2\dv \, \dd t-\int_0^T\!\!\!\int_M\varrho_2(x)\para{\sum_{j,k=1}^n \g^{jk}\frac{\p u_1}{\p x_j}
\frac{\p \overline{u}_2}{\p x_k}}\dv \, \dd t\cr
&\quad +\int_0^T\!\!\!\int_{\p\M}\p_\nu \overline{u}_2f_1\dsc \, \dd t
-\int_0^T\!\!\!\int_{\p\M}\p_\nu u_1\overline{f}_2\,\ds \dd t.
\end{align}
Taking into account the facts that $\Lambda_{c\g}$ is self-adjoint, that $c=1$ on $\p\M$ and $\para{\p_t^2\overline{u}_2-\Delta_{\cg}\overline{u}_2}=0$ in $(0,T)\times\M$,
it follows that
\begin{multline}\label{6.7}
 \int_{0}^T\!\!\!\int_{\p\M}\para{\Lambda_\g-\Lambda_{\cg}}f_1\,\overline{f}_2\ds \, \dd t = \hfill \cr
\int_0^T\!\!\!\int_\M\varrho_1(x)\p_t u_1\p_t\overline{u}_2\,\dv \,\dd t -\int_0^T\!\!\!\int_\M \varrho_2(x)\para{\sum_{j,k=1}^n \g^{jk}\frac{\p u_1}{\p x_j}\frac{\p \overline{u}_2}{\p x_k}}\dv \,\dd t.
\end{multline}
This completes the proof of the Lemma.
\end{Demo}
\subsection{Modified geometrical optics solutions}

As in the case of potentials, we extend the  manifold $(\M,\g)$ into a simple manifold
$\M_2 \Supset \M$ so that $\M_{2} \Supset \M_{1} \Supset \M$ with $(\M_1,\g)$ simple. We extend the conformal factor $c$
by $1$ outside the manifold $\M$; its $\mathcal{C}^k(\M_{1})$ norms  may also be bounded by $M_{0}$.
Let $\psi_1$, $\psi_2$ be two phase functions solving the eikonal equation with respect respectively to the metrics $\g$ and $\cg$.
\begin{align}\label{6.8}
\abs{\nabla_\g\psi_1}^2_\g=\sum_{j,k=1}^n\g^{jk}\frac{\p\psi_1}{\p x_j}\frac{\p\psi_1}{\p x_k}=1,\quad
\abs{\nabla_{\cg}\psi_2}^2_{\cg}=\sum_{j,k=1}^n{\cg}^{jk}\frac{\p\psi_2}{\p x_j}\frac{\p\psi_2}{\p x_k}=1.
\end{align}
Let $a_2$ solve the transport equation in $\R\times\M$ with respect the metric $\g$ (as given in section 4)
\begin{equation}\label{6.9}
\frac{\p a_2}{\p t}+\sum_{j,k=1}^n\g^{jk}\frac{\p\psi_1}{\p x_j}\frac{\p a_2}{\p x_k}+\frac{a_2}{2}\Delta_\g\psi_1=0.
\end{equation}
Let $a_3$ solve the following transport equation in $\R\times\M$ with respect to the metric $\cg$
\begin{align}\label{6.10}
\nonumber
\frac{\p a_3}{\p t}+\sum_{j,k=1}^n{\cg}^{jk}\frac{\p\psi_2}{\p x_j}\frac{\p a_3}{\p x_k}+\frac{a_3}{2}\Delta_{\cg}\psi_2&=-
\frac{1}{2i}a_2(t,x)(1-c^{-1})e^{i\lambda(\psi_1-\psi_2)}\\  &\equiv a_2(t,x)\varphi_0(x,\lambda),
\end{align}
and be such that
\begin{equation}\label{6.11}
\norm{a_3}_*\leq C\varepsilon\lambda^2\norm{a_2}_*.
\end{equation}
Let us explain the construction of a solution $a_3$ satisfying (\ref{6.10}) and (\ref{6.11}). To solve the transport equation (\ref{6.10})
and (\ref{6.11}) it is enough to take, in geodesic polar coordinates $(r,\theta)$ (with respect to the metric $\cg$)
\begin{equation}\label{6.12}
\widetilde{a}_3(t,r,\theta;\lambda)=\alpha_{\cg}^{-1/4}(r,\theta)\int_0^r\alpha_{\cg}^{1/4}(s,\theta)\widetilde{a}_2(s-r+t,s,\theta)
\widetilde{\varphi}_0(s,\theta,\lambda) \, \dd s,
\end{equation}
where $\alpha_{\cg}(r,\theta)$ denotes the square of the volume element in geodesic polar coordinates with respect to the metric $\cg$.
Using that $\norm{\varphi_0(\cdot,\lambda)}_{\mathcal{C}^2(\M)}\leq C\varepsilon\lambda^2$ and (\ref{6.12}) we obtain (\ref{6.11}).
\begin{Lemm}\label{L6.2}
Let $c\in \mathscr{C}(M_0,k,\varepsilon)$ be such that $c=1$ near the boundary $\p\M$. Then the equation
\begin{equation}\label{6.13}
\para{\p_t^2-\Delta_{\cg}}u=0,\quad\textrm{in}\quad (0,T)\times\M,\quad u(0,x)=\p_tu(0,x)=0
\end{equation}
has a solution of the form
\begin{equation}\label{6.14}
u_2(t,x)=\frac{1}{\lambda}a_2(t,x)e^{i\lambda(\psi_1(x)-t)}+a_3( t,x;\lambda)e^{i\lambda(\psi_2(x)- t)}+
v_{2,\lambda}(t,x)
\end{equation}
when $\lambda$ is large enough, which satisfies
\begin{multline}\label{6.15}
\lambda\norm{v_{2,\lambda}(t,\cdot)}_{L^2(\M)}+\norm{\nabla_{\g} v_{2,\lambda}(t,\cdot)}_{L^2(\M)}+\norm{\p_tv_{2,\lambda}(t,\cdot)}_{L^2(\M)}
\cr \leq C\para{\varepsilon\lambda^2+\lambda^{-1}}\norm{a_2}_*.
\end{multline}
The constant $C$ depends only on $T$ and $\M$ (that is $C$ does not depend on $a$, $\lambda$ and $\varepsilon$).
\end{Lemm}
\begin{Demo}{}
We set
\begin{align}\label{6.18}
k(t,x)=&-\para{\partial_t^2-\Delta_{\cg}}\para{\frac{1}{\lambda}a_2(t,x)e^{i\lambda(\psi_1- t)} + a_3(t,x,\lambda)e^{i\lambda(\psi_2- t)}}.
\end{align}
To prove our Lemma it would be enough to show that if $v$ solves
\begin{equation}\label{6.16}
\para{\partial_t^2-\Delta_{\cg}}v=k(t,x)
\end{equation}
with initial and boundary conditions
\begin{equation}\label{6.17}
v(0,x)=\p_tv(0,x)=0, \textrm{ in }\M,\quad \textrm{and}\quad v(t,x)=0 \textrm{ on }(0,T)\times\p\M
\end{equation}
then the estimates (\ref{6.15}) holds. We have
\begin{align}\label{6.19}
-k(t,x)&=\frac{1}{\lambda}e^{i\lambda(\psi_1-t)}\para{\p_t^2-\Delta_{\cg}}a_2\cr
& \quad
+2ie^{i\lambda(\psi_1-t)}\para{\p_ta_2+\sum_{j,k=1}^n
{\cg}^{jk}\frac{\p\psi_1}{\p x_j}\frac{\p a_2}{\p x_k}+\frac{a_2}{2}\Delta_{\cg}\psi_1}\cr
& \quad
+\lambda a_2e^{i\lambda(\psi_1-t)}\para{1-c^{-1}\sum_{j,k=1}^n\g^{jk}\frac{\p\psi_1}{\p x_j}\frac{\p\psi_1}{\p x_k}}\cr
& \quad
+e^{i\lambda(\psi_2-t)}\para{\p_t^2-\Delta_{\cg}}a_3\cr
& \quad
+2i\lambda e^{i\lambda(\psi_2- t)}\para{\p_ta_3+\sum_{j,k=1}^n{\cg}^{jk}\frac{\p\psi_2}{\p x_j}
\frac{\p a_3}{\p x_k}+\frac{a_3}{2}\Delta_{\cg}\psi_2}\cr
& \quad
+\lambda^2a_3 e^{i\lambda(\psi_2-t)}\para{1-\sum_{j,k=1}^n{\cg}^{jk}\frac{\p\psi_2}{\p x_j}\frac{\p\psi_2}{\p x_k}}.
\end{align}
Taking into account (\ref{6.8}) and (\ref{6.9}), the right-hand side of (\ref{6.19}) becomes
\begin{align}\label{6.20}
-k(t,x)&=\frac{1}{\lambda}e^{i\lambda(\psi_1-t)}\para{\p_t^2-\Delta_{\cg}} a_2 \cr
& \quad
+2ie^{i\lambda(\psi_1- t)}\Big((c^{-1}-1)\seq{\nabla_\g\psi_1,\nabla_\g a_2}_\g+
\frac{1}{2}a_2\para{\Delta_{\cg}\psi_1-\Delta_\g\psi_1}\Big)\cr
& \quad
\begin{aligned}
+2i\lambda e^{i\lambda(\psi_2- t)}\Big(\p_ta_3&+\sum_{j,k=1}^n
c{\g}^{jk}\frac{\p\psi_2}{\p x_j}\frac{\p a_3}{\p x_k}+\frac{a_3}{2}\Delta_{\cg}\psi_2 \\
&+\frac{a_2}{2i}e^{i\lambda(\psi_1-\psi_2)}(1-c^{-1})\Big)
\end{aligned} \cr
& \quad
+e^{i\lambda(\psi_2-t)}\para{\p_t^2-\Delta_{\cg}}a_3.
\end{align}
By (\ref{6.10}) we get
\begin{align}\label{6.21}
-k(t,x)&=\frac{1}{\lambda}e^{i\lambda(\psi_1-t)}\para{\p_t^2-\Delta_{\cg}}a_2\cr
& \quad
+2ie^{i\lambda(\psi_1- t)}\Big((c^{-1}-1)\seq{\nabla_\g\psi_1,\nabla_\g a_2}_\g  +
\frac{1}{2}a_2\para{\Delta_{\cg}\psi_1-\Delta_\g\psi_1}\Big) \cr
& \quad
+e^{i\lambda(\psi_2-t)}\para{\p_t^2-\Delta_{\cg}}a_3\cr
&
\equiv \frac{1}{\lambda}e^{i\lambda(\psi_1-t)}k_0+e^{i\lambda(\psi_1-\lambda t)}k_1
+e^{i\lambda(\psi_2-t)}k_2.
\end{align}
Since $k_j\in L^1(0,T;L^2(\M))$, by Lemma \ref{L1.1}, we deduce that
\begin{equation}\label{6.22}
v_\lambda\in {\cal C}^1(0,T;L^2(\M))\cap{\cal C}(0,T;
H^1_0(\M))
\end{equation}
and
\begin{align}\label{6.23}
\|&v_\lambda(t,\cdot)\|_{L^2(\M)} \cr
&\leq \frac{C}{\lambda} \bigg\{\int_{\R} \para{\frac{1}{\lambda}\norm{k_0(s,\cdot)}_{L^2(\M)}+\norm{k_1(s,\cdot)}_{L^2(\M)}+\norm{k_2(s,\cdot)}_{L^2(\M)}} \, \dd s\cr
& \quad
+\int_{\R} \para{\frac{1}{\lambda}\norm{\p_tk_0(s,\cdot)}_{L^2(\M)}+\norm{\p_tk_1(s,\cdot)}_{L^2(\M)}+\norm{\p_tk_2(s,\cdot)}_{L^2(\M)}} \, \dd s\bigg\}\cr
&
\leq \frac{C}{\lambda}\para{\frac{1}{\lambda}\norm{a_2}_*+\varepsilon\norm{a_2}_*+\varepsilon\lambda^2\norm{a_2}_*}\cr
&
\leq C\para{\varepsilon\lambda+\frac{1}{\lambda^2}}\norm{a_2}_*.
\end{align}
Moreover, we have
\begin{eqnarray}\label{6.24}
\norm{k}_{L^2((0,T)\times\M)}\leq C\para{\frac{1}{\lambda}\norm{a_2}_*+\varepsilon\norm{a_2}_*+\varepsilon\lambda^2\norm{a_2}_*}
\end{eqnarray}
and by using again the energy estimates for the problem (\ref{6.16})-(\ref{6.17}), we obtain
\begin{eqnarray}\label{6.25}
\norm{\nabla v_\lambda (t,\cdot)}_{L^2(\M)}+\norm{\p_t v_\lambda (t,\cdot)}_{L^2(\M)}\leq C\para{\varepsilon\lambda^2+\frac{1}{\lambda}}\norm{a_2}_*.
\end{eqnarray}
This ends the proof of Lemma \ref{L6.2}.
\end{Demo}
\begin{Lemm}\label{L6.3}
There exists $C>0$ such that for any $a_1$, $a_2\in H^1(\R,H^2(\M))$ satisfying the transport equation (\ref{6.9}) with (\ref{4.3}) the
following estimate holds true
\begin{multline}\label{6.28}
\abs{\int_0^T\!\!\!\int_\M \varrho(x)(a_1\overline{a}_2)(t,x)\dv dt}\leq \norm{\varrho_0}_{\mathcal{C}(\M)}
\para{\lambda^{-1}+\varepsilon\lambda^3}\norm{a_1}_*\norm{a_2}_*\cr
+\lambda^3\norm{a_1}_*\norm{a_2}_*\norm{\Lambda_{\g}-\Lambda_{\cg}}
\end{multline}
for any sufficiently large $\lambda$.
\end{Lemm}
\begin{Demo}{}
Following Lemma \ref{L6.2} let $u_2$ be a solution to the problem $(\p_t^2-\Delta_{\cg})u=0$ of the form
$$
\overline{u}_2(t,x)=\frac{1}{\lambda}\overline{a}_2( t,x)e^{-i\lambda(\psi_1- t)}+\overline{a}_3
( t,x;\lambda)e^{-i\lambda(\psi_2- t)}+\overline{v}_{2,\lambda}(t,x)
$$
where $v_{2,\lambda}$ satisfies (\ref{6.15}) and $a_3$ satisfies (\ref{6.11}).
Thanks to Lemma \ref{L4.1} let $u_1$ be a solution to the $(\p_t^2-\Delta_\g)u=0$ of the form
$$
u_1(t,x)=a_1( t,x)e^{i\lambda(\psi_1-t)}+v_{1,\lambda}(t,x),
$$
where $v_{1,\lambda}$ satisfies (\ref{4.7}). Then
\begin{align}\label{6.29}
\p_t\overline{u}_2(t,x)&= \frac{1}{\lambda}\p_t\overline{a}_2( t,x)e^{-i\lambda(\psi_1- t)}+i
\overline{a}_2( t,x)e^{-i\lambda(\psi_1- t)}\cr
& \quad +\p_t\overline{a}_3( t,x;\lambda)e^{-i\lambda(\psi_2- t)}+i\lambda\overline{a}_3
( t,x,\lambda)e^{-i\lambda(\psi_2- t)}
+\p_t\overline{v}_{2,\lambda}(t,x)\cr
\p_tu_1(t,x)&=\p_ta_1(t,x)e^{i\lambda(\psi_1-t)}-i\lambda a_1(t,x)e^{i\lambda(\psi_1-t)}+\p_tv_{1,\lambda}.
\end{align}
Let us compute the first term in the right hand side of (\ref{6.5}). We have
\begin{align}\label{6.30}
\nonumber &\int_0^T\!\!\!\int_\M \varrho_1\p_tu_1\p_t\overline{u}_2\dv \, \dd t=\lambda\int_0^T\!\!\!\int_\M \varrho_1a_1\overline{a}_2\dv \, \dd t \\ \nonumber
&+\frac{1}{\lambda}\int_0^T\!\!\!\int_\M \varrho_1\para{\p_ta_1\p_t\overline{a}_2}\dv \, \dd t
-i\int_0^T\!\!\!\int_\M \varrho_1 a_1\p_t\overline{a}_2 \dv \, \dd t \\ \displaybreak[1] \nonumber
&+\frac{1}{\lambda}\int_0^T\!\!\!\int_\M \varrho_1\p_t\overline{a}_2\p_tv_{1,\lambda} e^{-i\lambda(\psi_1-t)}\dv \, \dd t
+i\int_0^T\!\!\!\int_\M \varrho_1\para{\p_ta_1\overline{a}_2}\dv \, \dd t  \\ \nonumber
&+i \int_0^T\!\!\!\int_\M \varrho_1\overline{a}_2\p_tv_{1,\lambda} e^{-i\lambda(\psi_1-t)}\dv \, \dd t+\int_0^T\!\!\!\int_\M
\varrho_1\p_t\overline{a}_3\p_ta_{1}e^{i\lambda(\psi_1-\psi_2)}\dv \, \dd t  \\ \nonumber
&-i\lambda\int_0^T\!\!\!\int_\M \varrho_1\p_t\overline{a}_3a_{1}e^{i\lambda(\psi_1-\psi_2)}\dv \,\dd t+\!\int_0^T\!\!\!\int_\M
\varrho_1\p_tv_{1,\lambda}\p_t\overline{a}_3e^{-i\lambda(\psi_2-t)}\dv \dd t  \\ \nonumber
&+i\lambda\int_0^T\!\!\!\int_\M \varrho_1\overline{a}_3\p_ta_{1}e^{i\lambda(\psi_1-\psi_2)}\dv \, \dd t+\lambda^2\int_0^T\!\!\!\int_\M
\varrho_1a_1\overline{a}_3e^{i\lambda(\psi_1-\psi_2)}\dv \, \dd t  \\ \nonumber
&+i\lambda\int_0^T\!\!\!\int_\M \varrho_1\p_tv_{1,\lambda}\overline{a}_3e^{-i\lambda(\psi_2-t)}\dv \, \dd t+\!\int_0^T\!\!\!\int_\M
\varrho_1\p_ta_1\p_t\overline{v}_{2,\lambda}e^{i\lambda(\psi_1-t)} \dv \dd t  \\
&-i\lambda\int_0^T\!\!\!\int_\M \varrho_1a_1\p_t\overline{v}_{2,\lambda}e^{i\lambda(\psi_1-t)} \dv \, \dd t+\int_0^T\!\!\!\int_\M
\varrho_1\p_tv_{1,\lambda}\p_t\overline{v}_{2,\lambda}\dv \, \dd t.
\end{align}
Thus, we have from (\ref{6.11}), (\ref{6.15}) and (\ref{4.7}) the following identity
\begin{equation}\label{6.31}
\int_0^T\!\!\!\int_\M\varrho_1(x)\p_tu_1\p_t\overline{u}_2\dv dt=
\lambda\int_0^T\!\!\!\int_\M \varrho_1(x)(a_1\overline{a}_2)(t,x)\dv dt+\mathcal{J}_1(\lambda,\varepsilon),
\end{equation}
where
\begin{equation}\label{6.32}
\abs{\mathcal{J}_1(\lambda,\varepsilon)}\leq\norm{\varrho_0}_{\mathcal{C}(\M)}\para{1+\varepsilon\lambda^4}\norm{a_2}_*\norm{a_1}_*.
\end{equation}
On the other hand, we have
\begin{align}
\nonumber
\nabla_\g u_1 &= \nabla_\g a_1 e^{i\lambda(\psi_1-t)}+
i\lambda \nabla_\g \psi_1 a_1 e^{i\lambda(\psi_1- t)}+\nabla_\g v_{1,\lambda} \\
\label{6.33}
\nabla_\g \overline{u}_2 &=\frac{1}{\lambda}\nabla_\g\overline{a}_2 e^{-i\lambda(\psi_1- t)}-i\overline{a}_2 \nabla_\g\psi_1e^{-i\lambda(\psi_1- t)} \cr
&\quad -i\lambda\overline{a}_3\nabla_{\g}\psi_2e^{-i\lambda(\psi_2- t)}+\nabla_\g\overline{a}_3 e^{-i\lambda(\psi_2- t)}+\nabla_\g\overline{v}_{2,\lambda}.
\end{align}
and the second term in the right side of (\ref{6.5}) becomes
\begin{multline}\label{6.35}
\int_0^T\!\!\!\int_\M \varrho_2(x)\seq{\nabla_\g u_1(t,x),\nabla_\g \overline{u}_2(t,x)}_\g \dv dt
\\ =\lambda\int_0^T\!\!\!\int_\M \varrho_2(x)(a_1\overline{a}_2)(t,x)\dv dt+\mathcal{J}_2(\lambda,\varepsilon)+\mathcal{J}_3(\lambda,\varepsilon)
\end{multline}
with
\begin{align*}
\mathcal{J}_2(\lambda,\varepsilon)=
&+\frac{1}{\lambda}\int_0^T\!\!\!\int_\M\varrho_2\seq{\nabla_\g a_1,\nabla_\g\overline{a}_2}_\g\dv \, \dd t\\
&-i\int_0^T\!\!\!\int_\M\varrho_2\overline{a}_2\seq{\nabla_\g a_1,\nabla_\g\psi_1(x)}_\g\dv \, \dd t\\
&+i\int_0^T\!\!\!\int_\M \varrho_2 a_1\seq{\nabla_\g\overline{a}_2, \nabla_\g\psi_1}_\g\dv \dd t\\  \displaybreak[1]
&+\frac{1}{\lambda}\int_0^T\!\!\!\int_\M \varrho_2 e^{-i\lambda(\psi_1- t)}\seq{\nabla_\g\overline{a}_2,\nabla_\g v_{1,\lambda}}_\g\dv \, \dd t\\
&-i\int_0^T\!\!\!\int_\M \varrho_2 \overline{a}_2 e^{-i\lambda(\psi_1- t)}\seq{\nabla_\g v_{1,\lambda},\nabla_\g\psi_1}_\g\dv \, \dd t
\end{align*}
and
\begin{align*}
\mathcal{J}_3(\lambda,\varepsilon)=
&-i\lambda\int_0^T\!\!\!\int_\M \varrho_2\overline{a}_3 e^{i\lambda(\psi_1-\psi_2)}\seq{\nabla_\g a_1,\nabla_\g\psi_2}_\g\dv \, \dd t  \displaybreak[1]\\
&+\int_0^T\!\!\!\int_\M\varrho_2 e^{i\lambda(\psi_1-\psi_2)}\seq{\nabla_\g a_1,\nabla_\g\overline{a}_3}_\g\dv \, \dd t  \displaybreak[1] \\
&+\int_0^T\!\!\!\int_\M \varrho_2 e^{i\lambda(\psi_1- t)}\seq{\nabla_\g a_1,\nabla_\g\overline{v}_{2,\lambda}}_\g\dv \, \dd t  \displaybreak[1] \\
&+\lambda^2\int_0^T\!\!\!\int_\M \varrho_2 a_1\overline{a}_3 e^{i\lambda(\psi_1-\psi_2)}\seq{\nabla_\g\psi_1,\nabla_\g\psi_2}_\g\dv \, \dd t  \displaybreak[1] \\
&+i\lambda\int_0^T\!\!\!\int_\M \varrho_2 a_1e^{i\lambda(\psi_1-\psi_2)}\seq{\nabla_\g\overline{a}_3,\nabla_\g\psi_1}_\g\dv \, \dd t \displaybreak[1] \\
&+i\lambda\int_0^T\!\!\!\int_\M \varrho_2 a_1e^{i\lambda(\psi_1- t)}\seq{\nabla_\g\psi_1,\nabla_\g\overline{v}_2}_\g \dv \, \dd t  \displaybreak[1] \\
&-i\lambda\int_0^T\!\!\!\int_\M \varrho_2\overline{a}_3e^{-i\lambda(\psi_2- t)}\seq{\nabla_\g v_{1,\lambda},\nabla_g\psi_2}_\g\dv \, \dd t  \displaybreak[1] \\
&+\int_0^T\!\!\!\int_\M \varrho_2 e^{-i\lambda(\psi_2- t)}\seq{\nabla_\g v_{1,\lambda},\nabla_\g \overline{a}_3}_\g\dv \, \dd t  \displaybreak[1] \\
&+\int_0^T\!\!\!\int_\M \varrho_2 \seq{\nabla_\g v_{1,\lambda}, \nabla_\g \overline{v}_{2,\lambda}}_\g\dv \, \dd t.
\end{align*}
From (\ref{6.11}), (\ref{6.15}) and (\ref{4.7}), we have
\begin{equation}\label{6.36}
\abs{\mathcal{J}_2(\lambda,\varepsilon)}+\abs{\mathcal{J}_3(\lambda,\varepsilon)}
\leq \norm{\varrho_0}_{\mathcal{C}(\M)}\para{1+\varepsilon\lambda^4}\norm{a_2}_*\norm{a_1}_*.
\end{equation}
Taking into account (\ref{6.5}), (\ref{6.31}) and (\ref{6.35}), we deduce that
\begin{multline}\label{6.37}
\int_0^T\!\!\!\int_{\p\M}\para{\Lambda_{\g}-\Lambda_{\cg}}f_1\overline{f}_2\ds dt \\ =
\lambda\int_0^T\!\!\!\int_\M \varrho(x)(a_1\overline{a}_2)(t,x)\dv dt+\mathcal{J}_1(\lambda,\varepsilon)+
\mathcal{J}_2(\lambda,\varepsilon)+\mathcal{J}_3(\lambda,\varepsilon).
\end{multline}
In view of (\ref{6.36}) and (\ref{6.32}), we obtain
\begin{multline}\label{6.38}
\abs{\int_0^T\!\!\!\int_\M \varrho(x)(a_1\overline{a}_2)( t,x)\dv dt}\leq \norm{\varrho_0}_{\mathcal{C}(\M)}
\para{\lambda^{-1}+\varepsilon\lambda^3}\norm{a_1}\norm{a_2}\\
+\lambda^3\norm{a_1}_*\norm{a_2}_*\norm{\Lambda_{\g}-\Lambda_{\cg}}.
\end{multline}
This completes the proof.
\end{Demo}
\subsection{Stability estimate of the geodesic ray transform}
\begin{Lemm}\label{L6.4}
Let $M_0>0$. There exist $C>0$ and $\beta_j>0$, $j=1,2,3$, such that for any $b\in H^2(\p_+S\M_{1})$ the following estimate
\begin{multline}\label{6.39}
\abs{\int_{\p_+S\M_1}\I(\varrho)(y,\theta)b(y,\theta)\seq{\theta,\nu(y)}\dss (y,\theta)}\cr
\leq \Big((\lambda^{-\beta_1}+\varepsilon\lambda^{\beta_2})\norm{\varrho_0}_{\mathcal{C}^1(\M)}+\lambda^{\beta_3}
\norm{\Lambda_{\g}-\Lambda_{\cg}}\Big)\norm{b}_{H^2(\p_+S\M_{1})}.
\end{multline}
holds for any $\lambda$ large.
\end{Lemm}
\begin{Demo}{}
Following (\ref{4.21}), we take two solutions of the form
\begin{align*}
\widetilde{a}_1(t,r,\theta)&=\alpha^{-1/4}\phi(t-r)b(y,\theta), \\
\widetilde{a}_2(t,r,\theta)&=\alpha^{-1/4}\phi(t-r)\mu(y,\theta).
\end{align*}
Now we change variable in (\ref{6.28}), $x=\exp_{y}(r\theta)$, $r>0$ and
$\theta\in S_{y}\M_1$. Then
\begin{align*}
\int_0^T\!\!\int_\M &\varrho(x)a_1( t,x)a_2( t,x)
\dv\,dt\cr
&=\int_0^T\!\!\int_{S_{y}\M_1}\!\int_0^{\tau_+(y,\theta)}\widetilde{\varrho}(r,\theta)\widetilde{a}_1( t,r,\theta)
\widetilde{a}_2( t,r,\theta)
\alpha^{1/2} \,\dd r \,\dd\omega_y(\theta)\,\dd t \cr
&=\int_0^T\!\!\int_{S_{y}\M_1}\!\int_0^{\tau_+(y,\theta)}\widetilde{\varrho}(r,\theta)\phi^2( t-r)b(y,\theta)\mu(y,\theta) \,\dd r \,\dd\omega_y(\theta)\,\dd t.
\end{align*}
We conclude that
\begin{multline}\label{6.41}
\abs{\int_0^{\infty}\!\!\!\int_{S_{y}\M_1}\!\int_0^{\tau_+(y,\theta)}\widetilde{\varrho}(r,\theta)\phi^2( t-r)b(y,\theta)\mu(y,\theta)
 \,\dd r \,\dd\omega_y(\theta)\,\dd t}\\ \leq C\Big((\lambda^{-1}+\varepsilon\lambda^3)\norm{\varrho_0}_{\mathcal{C}(\M)}+\lambda^3
\norm{\Lambda_{\g,\,q_1}-\Lambda_{\g,\,q_2}}\Big) \\ \times \norm{\phi}^2_{H^3(\R)}\norm{b(y,\cdot)}_{H^2(S_y^+\M_{1})}.
\end{multline}
Since $\phi(t)=0$ for $t\leq 0$ or $t\geq T$, we get
\begin{multline}\label{6.41'}
\int_0^{\infty}\!\!\!\int_{S_{y}\M_1}\!\int_0^{\tau_+(y,\theta)}\widetilde{\varrho}(r,\theta)\phi^2( t-r)b(y,\theta)\mu(y,\theta) \, \dd r \, \dd \omega_y(\theta) \, \dd t \\ =
\bigg(\int_{-\infty}^{\infty} \phi^2(t) dt\bigg) \times \int_{S_{y}\M_1}\!\int_0^{\tau_+(y,\theta)}\widetilde{\varrho}(r,\theta)b(y,\theta)\mu(y,\theta) \, \dd r \, \dd \omega_y(\theta).
\end{multline}
Combining (\ref{6.41'}) and (\ref{6.41}), it follows that
\begin{multline*}
\abs{\int_{S_{y}\M_1}\!\int_0^{\tau_+(y,\theta)}\widetilde{\varrho}(s,\theta)b(y,\theta)\mu(y,\theta) \, \dd s \, \dd\omega_y(\theta) } \\ \leq
 C\Big((\lambda^{-1}+\varepsilon\lambda^3)\norm{\varrho_0}_{\mathcal{C}(\M)}+\lambda^3
\norm{\Lambda_{\g,\,q_1}-\Lambda_{\g,\,q_2}}\Big)\norm{b(y,\cdot)}_{H^2(S_y^+\M_{1})}.
\end{multline*}
Integrating with respect to $y\in \p \M_1$ we obtain
\begin{multline}\label{6.46}
\abs{\int_{\p_+S\M_1}\I(\varrho)(y,\theta)b(y,\theta)\seq{\theta,\nu(y)}\dss (y,\theta)}\\
\leq \Big((\lambda^{-\beta_1}+\varepsilon\lambda^{\beta_2})\norm{\varrho_0}_{\mathcal{C}^1(\M)}+\lambda^{\beta_3}
\norm{\Lambda_{\g}-\Lambda_{\cg}}\Big)\norm{b}_{H^2(\p_+S\M_{1})}.
\end{multline}
This completes the proof of the lemma.
\end{Demo}
\subsection{End of the proof of Theorem \ref{Th2}}
This subsection is devoted to the end of the proof of Theorem \ref{Th2}. We need the following known result (see \cite{[SU2]} proposition 4.1).
\begin{Lemm}\label{L6.5}
Let $c\in\mathcal{C}^\infty(\M)$ be such that $\norm{1-c}_{\mathcal{C}(\M)}\leq\varepsilon$. Then there exists $C>0$ such that
\begin{equation}\label{6.47}
\norm{d_\g-d_{\cg}}_{\mathcal{C}(\p\M\times\p\M)}\leq C\norm{\Lambda_\g-\Lambda_{\cg}}^\mu,
\end{equation}
with some $0<\mu<1$ depending only on the dimension $n$.
\end{Lemm}
From this lemma, we can derive the following estimate
\begin{Corol}\label{C6.1} There exists a constant $C>0$ such that the following estimate holds true
\begin{equation}\label{6.48}
\norm{\I^*\I(\varrho)}_{H^2(\M_1)}^{3/2}\leq C\para{\norm{\varrho}_{\mathcal{C}^1(\M)}^2+\norm{\Lambda_\g-\Lambda_{\cg}}^\mu}
\end{equation}
with some $0<\mu<1$ depending on $n$ only.
\end{Corol}
\begin{Demo}{}
Linearizing near $\g$, we get, as in \cite{[Eskin4]}
$$
d_\g(x,y)-d_{\cg}(x,y)=\frac{1}{2}\I(\varrho)(x,y)+\mathcal{R}(\varrho)(x,y),\quad \forall x,y\in\p\M,
$$
where, with some abuse of notation, $\I(\varrho)(x,y)$ stands for $\I(\varrho)(x,\theta)$ with $\theta=\exp^{-1}_x(y)/\abs{\exp^{-1}_x(y)}$.
The remainder term $\mathcal{R}(\varrho)(x,y) $ is nonlinear and satisfies the estimate (see \cite{[Eskin4]})
$$
\abs{\mathcal{R}(\varrho)(x,y)}\leq Cd_{\g}(x,y) \norm{\varrho}_{\mathcal{C}^1(\M)}^2,\quad \forall x,y\in\p\M,
$$
with $C>0$ uniform in $c$ if $0<\varepsilon\ll 1$. By Lemma \ref{L6.5}, we have
$$
\abs{\I(\varrho)(x,y)}\leq C\para{d_{\g}(x,y) \norm{\varrho}_{\mathcal{C}^1(\M)}^2+\norm{\Lambda_\g-\Lambda_{\cg}}^\mu},\quad \forall x,y\in\p\M.
$$
Apply $\I^*$ to both sides, and use the estimate $\norm{\I^*(f)}_{L^\infty(\M_1)}\leq C\norm{f}_{L^\infty(\M_1)}$ to get
\begin{equation}\label{6.48'}
\norm{\I^*\I(\varrho)}_{L^\infty(\M_1)}\leq C\para{\norm{\varrho}_{\mathcal{C}^1(\M)}^2+\norm{\Lambda_\g-\Lambda_{\cg}}^\mu}.
\end{equation}
Since $\varrho$ vanishes outside $\M$ with all derivatives and $\I^*\I$ is a pseudodifferential operator of order $-1$, we have
$$
\norm{\I^*\I(\varrho)}_{H^{m+1}(\M_1)}\leq C_m\norm{\varrho}_{H^m(\M)}
$$
for all integers $m$. Using interpolation, we get
\begin{align}\label{6.48"}
\norm{\I^*\I(\varrho)}_{H^{2}(\M_1)}&\leq C\norm{\I^*\I(\varrho)}_{L^2(\M_1)}^{2/3}\norm{\I^*\I(\varrho)}_{H^{6}(\M_1)}^{1/3} \cr
&\leq C\norm{\I^*\I(\varrho)}_{L^\infty(\M_1)}^{2/3}.
\end{align}
Therefore, (\ref{6.48"}) and (\ref{6.48'}) imply
$$
\norm{\I^*\I(\varrho)}_{H^{2}(\M_1)}^{3/2}\leq C\para{\norm{\varrho}_{\mathcal{C}^1(\M)}^2+\norm{\Lambda_\g-\Lambda_{\cg}}^\mu}.
$$
This completes the proof.
\end{Demo}
Let us now prove Theorem \ref{Th2}. We choose
$$
b(y,\theta)=\I\para{\I^*\I(q)}(y,\theta)
$$
and obtain
\begin{multline*}
\norm{\I^*\I(\varrho)}^2_{L^2(\M_1)} \\ \leq C\Big((\lambda^{-\beta_1}+\varepsilon\lambda^{\beta_2})\norm{\varrho_0}_{\mathcal{C}^1(\M)}+\lambda^{\beta_3}
\norm{\Lambda_{\g}-\Lambda_{\cg}}\Big)\norm{\I^*\I(\varrho)}_{H^2(\M_1)}.
\end{multline*}
By interpolation we have
\begin{align}\label{6.49}
\|\I^*\I(&\varrho)\|^2_{H^1(\M_1)} \cr &\leq C \norm{\I^*\I(\varrho)}_{L^2(\M_1)}\norm{\I^*\I(\varrho)}_{H^2(\M_1)}\cr \nonumber
&\leq C \Big((\lambda^{-\beta_1}+\varepsilon\lambda^{\beta_2})\norm{\varrho_0}_{\mathcal{C}^1(\M)}+\lambda^{\beta_3}
\norm{\Lambda_{\g}-\Lambda_{\cg}}\Big)^{1/2} \norm{\I^*\I(\varrho)}_{H^2(\M_1)}^{3/2}.
\end{align}
Using (\ref{6.48}) we obtain
$$
\norm{\varrho}^2_{L^2(\M)}\leq C\Big((\lambda^{-\beta'_1}+\varepsilon\lambda^{\beta'_2})\norm{\varrho_0}^{5/2}_{\mathcal{C}^1(\M)}+\lambda^{\beta'_3}
\norm{\Lambda_{\g}-\Lambda_{\cg}}^\mu\Big).
$$
Since
$$
\norm{\varrho_0}_{\mathcal{C}^1(\M)}\leq C\norm{\varrho_0}_{H^{n/2+1+\epsilon}(\M)}\leq C\norm{\varrho_0}_{L^2(\M)}^{4/5}
\norm{\varrho_0}_{H^s(\M)}^{1/5}\leq C\norm{\varrho_0}_{L^2(\M)}^{4/5}
$$
we obtain
$$
\norm{\varrho}^2_{L^2(\M)}\leq C\Big((\lambda^{-\beta'_1}+\varepsilon\lambda^{\beta'_2})\norm{\varrho_0}^2_{L^2(\M)}+\lambda^{\beta'_3}
\norm{\Lambda_{\g}-\Lambda_{\cg}}^\mu\Big).
$$
Minimising $(\lambda^{-\beta'_1}+\varepsilon\lambda^{\beta'_2})$ with respect to $\lambda>0$, we get
$$
C'\norm{\varrho_0}^2_{L^2(\M)}\leq \norm{\varrho}^2_{L^2(\M)}\leq C\Big(\varepsilon^{\gamma}\norm{\varrho_0}^2_{L^2(\M)}+C_\varepsilon
\norm{\Lambda_{\g}-\Lambda_{\cg}}^\mu\Big).
$$
for $\varepsilon>0$ small enough we conclude and obtain (\ref{1.14}).

\section{Proof of Theorem \ref{Th3}}
\setcounter{equation}{0}

Let $q\in L^\infty(\M)$, we first define the elliptic Dirichlet-to-Neumann map; let $\sigma (A_q)=\{ \lambda _{k,q}\}$ be the spectrum of $A_q$ and $\rho
 (A_q)=\mathbb{C}\setminus  \sigma (A_q)$ be the resolvent set of
 $A_q$. From well known results (e.g., \cite{[Lions-Magenes]}), for any $z \in \rho (A_q)$
 and $h\in H^{3/2}(\p\M)$, the nonhomogeneous boundary value problem
$$
\left\{
\begin{array}{ll}
(-\Delta_\g  +q)u =z u,\quad &\mbox{in}\; \M
\\
u=h, &\mbox{on}\; \p\M
\end{array}
\right.
$$
has an unique solution in $u_{q,h}\in H^2(\M )$ and the Dirichlet-to-Neumann map
$$
\Pi_{\g,q}(z) : f\rightarrow \partial _\nu u_{q,h}|_{\p\M}
$$
defines a bounded operator from $H^{3/2}(\p\M )$ to
$H^{1/2}(\p\M)$. We fix $T> \mathrm{Diam}_\g\M$ and consider the following function space
$$
\mathcal{H}_1 =\set{ f\in H^{2n+4}(0,T;H^{\frac{3}{2}} (\p\M ));\;
 \partial _t ^jf(0,\cdot )=0,\; 0\leq j\leq 2n+3 \;},
$$
and the operator
$$
\mathscr{R}_{\g,q} f=\sum_{k\geq 1}\frac{1}{\lambda
_{q,k}^{n+2}}(\partial_\nu \phi
_{q,k})\int_0^t\frac{\sin \sqrt{\lambda
_{q,k}}(t-s)}{\sqrt{\lambda_{q,k}}}\langle
 -\partial _s^{2(n+2)}f(\cdot ,s), \partial _\nu \phi
 _{q,k}\rangle \, \dd s,
$$
where $\seq{ \cdot,\cdot}$ denotes the $L^2(\p\M)$-scalar product.
Then $\mathscr{R}_{\g,q}$ defines a bounded operator from
$\mathcal{H}_1$ to $\mathcal{H}_2=L^2(0,T;H^s(\p\M))$.
\medskip

We will need in the sequel the following three lemmas. Their proof
can be found in \cite{[AS]}  or can be deduced easily from the results
in this reference (see also \cite{[Ch]}). We fix $0\leq s<\frac{1}{2}$.

\begin{Lemm}\label{LII.2.1}
Let $q\in L^\infty (\M)$. Then for any $m>\frac{n}{2}$, $h\in
H^{3/2}(\p\M )$ and $z \in \rho (A_q)$, we have
$$
\frac{\p^m}{\p z^m}\Pi_{\g,q}(z )h=-m!\sum_{k\geq 1}
\frac{1}{(\lambda _{k,q}-z )^{m+1}}\langle h,\partial _\nu
 \phi _{k,q}\rangle \partial _\nu \phi _{k,q} .
$$
\end{Lemm}
\begin{Lemm}\label{LII.2.2}
Let $N$ be a non negative integer and let $q_1$, $q_2\in L^\infty
(\M)$ satisfy $0\leq q_1,q_2\leq M_0$ for some positive
constant $M$. Then there exists a positive constant $C$, depending
only on $\M$ and $M_0$, such that
$$
\norm{\frac{\p^j}{\p z^j}\Big[\Pi_{\g,q_1}(z
)-\Pi_{\g,q_2}(z)\Big]} _s\leq \frac{C}{\abs{z}^{p+\frac{1-2s}{4}}},\quad z \leq
 0\quad\mbox{and}\quad 0\leq j\leq  N,
$$
where $\norm{\,\cdot\,}_s$ denotes the norm in  ${\cal L}(H^{3/2}(\p\M);\, H^s(\p\M))$.
\end{Lemm}
\begin{Lemm}\label{LII.2.3}
For each  $f\in \mathcal{H}_1$, we have
\begin{equation}\label{II.2.1}
\Lambda^\sharp_{\g,q}f=\sum_{j=0}^{n+1}
\cro{\frac{\p^j}{\p z ^j}\Pi_{\g,q}(z)}_{|z=0}\para{-\partial _t^2f}+\mathscr{R}_{\g,q}f,
\end{equation}
where $\Lambda^\sharp_{\g,q}$ is the restriction of $\Lambda_{\g,q}$ to $\mathcal{H}_1$.
\end{Lemm}
First, we remark that  for
$q\in L^\infty(\M)$ and $f\in \mathcal{H}_1$ the problem (\ref{1.2}) has an unique solution
\begin{equation}
u\in L^2(0,T;H^2(\M))\cap H^2\para{0,T;L^2(\M)}.
\end{equation}
Moreover $\Lambda^\sharp_{\g,q}$ is a linear and continuous map from $\mathcal{H}_1$ into $\mathcal{H}_2$. Indeed, for $f\in \mathcal{H}_1$,
let $v$ solve the problem
\begin{equation}\label{a}
\left\{
\begin{array}{llll}
\para{\partial^2_t-\Delta_\g}v =0,  & \textrm{in }\; (0,T)\times \M,\cr
\\
v(0,\cdot )=0,\quad\p_tv(0,\cdot)=0 & \textrm{in }\; \M,\cr
\\
v=f, & \textrm{on} \,\,(0,T)\times\p \M.
\end{array}
\right.
\end{equation}
Then
$$
v\in L^2(0,T;H^2(\M))\cap H^2\para{0,T;L^2(\M)}.
$$
Furthermore
\begin{equation}\label{b}
\norm{v}_{L^2(0,T;H^2(\M))}\leq C\norm{f}_{\mathcal{H}_1}.
\end{equation}
Estimate (\ref{b}) is essentially known, but we give the proof for the readers' convenience. Let $v_1=\p_t^2v$. Then, by hyperbolic estimates, we have
$$
\norm{v_1}_{L^2(0,T;L^2(\M))}\leq C\norm{f}_{\mathcal{H}_1}.
$$
On the other hand, since $\Delta_\g v=v_1$, by the elliptic regularity, we get
$$
\norm{v}_{L^2(0,T;H^2(\M))}\leq C\para{\norm{v_1}_{L^2(0,T;L^2(\M))}+\norm{f}_{\mathcal{H}_1}}.
$$
Thus, we get
$$
\norm{\frac{\p v}{\p\nu}}_{\mathcal{H}_2}\leq C\norm{f}_{\mathcal{H}_1}.
$$
Now, for $q\in L^\infty(\M)$, let $w$ solve
\begin{equation}\label{12}
\left\{
\begin{array}{llll}
\para{\partial^2_t-\Delta_\g +q(x)}w=-q(x)v,  & \textrm{in }\; (0,T)\times \M,\cr
\\
w(0,\cdot )=0,\quad\p_tw(0,\cdot)=0 & \textrm{in }\; \M,\cr
\\
w=0, & \textrm{on} \,\,(0,T)\times\p \M,
\end{array}
\right.
\end{equation}
we apply Lemma \ref{L1.1} to $\p_tw$, we can prove
$$
\norm{w}_{L^2(0,T;H^2(\M))}\leq C\norm{qv}_{H^1(0,T;L^2(\M))}\leq C\norm{f}_{\mathcal{H}_1}.
$$
Thus, for $u=v+w$, we have
\begin{equation}\label{c}
\left\{
\begin{array}{llll}
\para{\partial^2_t-\Delta_\g +q(x)}u=0,  & \textrm{in }\; (0,T)\times \M,\cr
\\
u(0,\cdot )=0,\quad\p_tu(0,\cdot)=0 & \textrm{in }\; \M,\cr
\\
u=f, & \textrm{on} \,\,(0,T)\times\p \M,
\end{array}
\right.
\end{equation}
where
$$
u\in L^2(0,T;H^2(\M))\cap H^2(0,T;L^2(\M)),
$$
and
\begin{equation}\label{d}
\|\Lambda_{\g,q}^\sharp(f)\|_{\mathcal{H}_2}\leq C\norm{f}_{\mathcal{H}_1}.
\end{equation}
We shall denote by $\|\Lambda_{\g,q}^\sharp\|_{\mathcal{L}\para{\mathcal{H}_1,\mathcal{H}_2}}$ the operator norm of $\Lambda_{\g,q}^\sharp$.
\begin{Lemm}
Let $(\M,\g)$ be a simple Riemannian compact manifold with boundary of dimension $n \geq 2$, let $T> \mathrm{Diam}_\g(\M)$, there exist constants
$C > 0$ and $\kappa\in (0,1)$ such that for any real valued potentials $q_1,\,q_2\in\mathscr{Q}(M_0)$ such that $q_{1}=q_{2}$ on the boundary $\p \M$,
we have
\begin{equation}\label{e}
\norm{q_1-q_2}_{L^2(\M)}\leq C
\|\Lambda^\sharp_{\g,q_1}-\Lambda^\sharp_{\g,q_2}\|_{\mathcal{L}\para{\mathcal{H}_1,\mathcal{H}_2}}^{\kappa}
\end{equation}
where $C$ depends on $\M$, $T$, $M_0$, $n$, and $s$.
\end{Lemm}
\begin{Demo}{}
As in (\ref{5.11}), we have
\begin{eqnarray}\label{5.11'}
&&\abs{\int_0^T\!\!\!\int_{\p \M}\para{\Lambda^\sharp_{\g,\,q_1}-\Lambda^\sharp_{\g,\, q_2}}(f_{\lambda})\overline{g}_\lambda \ds\, dt}\cr
&&\qquad\leq \|\Lambda^\sharp_{\g,\,q_1}-\Lambda^\sharp_{\g,\, q_2}\|_{\mathcal{L}\para{\mathcal{H}_1,\mathcal{H}_2}}
\norm{f_\lambda}_{\mathcal{H}_1}\norm{g_\lambda}_{L^2((0,T)\times\p \M)}\cr
&&\qquad\leq  C\lambda^{2n+5}\norm{a_1}_{**}\norm{a_2}_{**}\|\Lambda^\sharp_{\g,\,q_1}-\Lambda^\sharp_{\g,\,
q_2}\|_{\mathcal{L}\para{\mathcal{H}_1,\mathcal{H}_2}}.
\end{eqnarray}
Where
$$
\norm{a}_{**}=\norm{a}_{H^{2n+4}(0,T;H^2(\M))}.
$$
Thus, we can complete the proof of (\ref{e}) in the same way as in section 5.2.
\end{Demo}

\medskip

We set $P(z)=\para{\Pi _{\g,q_1}(z)-\Pi_{\g,q_2}(z)}$, from Taylor's formula, we deduce for $1\leq j\leq n$, and $z \leq 0$
\begin{equation}\label{II.2.9}
P^{(j)}(0)=\sum_{p=j}^n\frac{(-z)^{p-j}}{(p-j)!}P^{(p)}(z)+\int_z^0\frac{(-\tau )^{n-j}}{(n-j)!}P^{(n+1)}(\tau ) \, \dd \tau .
\end{equation}
\begin{Lemm}\label{LII.2.4} There exist $C>0$ and $\mu_1\in(0,1)$ such that the following estimate
\begin{equation}\label{II.2.10}
\norm{ P^{(n+1)}(z) }_s\leq C\epsilon^{\mu_1}
\end{equation}
holds true for any $z\leq 0$. Here $C$ is a positive constant
depending on $M_0$, $\M$ and $\norm{\,\cdot\,}_s$
denotes the norm in $\mathcal{L}(H^{3/2}(\p\M);H^s(\p\M))$.
\end{Lemm}
\begin{Demo}{}
Let $h\in H^{3/2}(\p\M)$. It follows from Lemma \ref{LII.2.1}
\begin{eqnarray*}
P^{(n+1)}(z)h&=&-(n+1)!\sum_{k\geq 1} \frac{1}{(\lambda
_{k,q_1}-z)^{n+2}}\langle h,\partial _\nu \phi
_{k,q_1}\rangle \partial _\nu \phi _{k,q_1} \\
&& +(n+1)!\sum_{k\geq 1} \frac{1}{(\lambda _{k,q_2}-z
)^{n+2}}\langle h,\partial_\nu \phi _{k,q_2} \rangle \partial
_\nu \phi _{k,q_2}.
\end{eqnarray*}
We split $P^{(n+1)}(z)h$ into three terms $P^{(n+1)}(z)h=\mathcal{I}_1(z)h+\mathcal{I}_2(z)h+\mathcal{I}_3(z)h$, where
\begin{align*}
\mathcal{I}_1(z)h &= -(n+1)!\sum_{k\geq 1}
\Big[\frac{1}{(\lambda _{k,q_1}-z)^{n+2}} - \frac{1}{(\lambda
 _{k,q_2}-z )^{n+2}}\Big]\langle h,\partial _\nu \phi _{k,q_1}
\rangle \partial _\nu \phi _{k,q_1}
\\
\mathcal{I}_2(z)h &= -(n+1)!\sum_{k\geq 1} \frac{1}{(\lambda
_{k,q_2}-z)^{n+2}}\langle h,\partial _\nu
 \phi _{k,q_1}
-\partial _\nu \phi _{k,q_2}\rangle \partial _\nu \phi
 _{k,q_1}
\\
\mathcal{I}_3(z)h &= -(n+1)!\sum_{k\geq 1} \frac{1}{(\lambda
_{k,q_2}-z)^{n+2}}\langle h,\partial _\nu \phi
_{k,q_2}\rangle [\partial _\nu \phi _{k,q_1}-\partial _\nu
\phi _{k,q_2}].
\end{align*}
For $\mathcal{I}_1(z)h$, we have
\begin{multline}
\norm{\mathcal{I}_1(z)h}_{H^{1/2}(\p\M)} \leq  (n+1)!\norm{h}_{L^2(\p\M)} \\ \times \sum_{k\ge 1}\abs{\frac{1}{(\lambda _{k,q_1}-z)^{n+2}} -
\frac{1}{(\lambda_{k,q_2}-z)^{n+2}}}\norm{\partial _\nu \phi _{k,q_2}}_{H^{1/2}(\p\M)}^2.
\end{multline}
On the other hand, noting  that $z \le 0$, $\lambda_{k,q_j}\ge 0$,
$j=1,2$, we see that
\begin{align*}
\abs{\frac{1}{(\lambda _{k,q_1}-z)^{n+2}} -
\frac{1}{(\lambda_{k,q_2}-z)^{n+2}}}&\leq C\max
\Big(\frac{1}{\lambda _{k,q_1}^{n+3}} , \frac{1}{\lambda_{k,q_2}^{n+3}}\Big) \abs{\lambda _{k,q_1}-\lambda _{k,q_2}}
\\
&\leq \frac{C}{k^{ \frac{2(n+3)}{n}} }\abs{\lambda_{k,q_1}-\lambda_{k,q_2}},
\end{align*}
where we have used estimate (\ref{II.1.2}). On the other hand, since (see (\ref{II.1.1}) and (\ref{II.1.2}))
$$
\norm{ \partial _\nu \phi _{k,q_2}}^2 _{H^{1/2}(\p\M)}\leq Ck^{\frac{4}{n}},
$$
we obtain
\begin{multline}\label{II.2.11}
\norm{\mathcal{I}_1(z)h}_{H^{1/2}(\p\M)} \\ \leq  C\norm{h}_{L^2(\p\M)}\sum_{k\geq 1}\frac{1}{k^{
 \frac{2(n+1)}{n}} }\abs{\lambda_{k,q_1}-\lambda_{k,q_2}}\leq C\epsilon\norm{h}_{L^2(\p\M)}.
\end{multline}
For $\mathcal{I}_2(z)h$, we have
\begin{multline}\label{II.2.12}
\norm{\mathcal{I}_2(z)h}_{H^{1/2}(\p\M)} \\ \leq
C\norm{h}_{L^2(\p\M)}\sum_{k\ge
1}\frac{\lambda_{k,q_1}}{(\lambda_{k,q_2}-z)^{n+2}}\norm{\p_\nu(\phi_{k,q_1}-\phi_{k,q_2})}_{L^2(\p\M)}.
\end{multline}
Then Lemma \ref{LII.2.4} yields
\begin{multline}\label{II.2.13}
\sum_{k\ge 1}\frac{\lambda_{k,q_1}}{(\lambda_{k,q_2}-z)^{n+2}}\norm{\p_\nu(\phi_{k,q_1}-\phi_{k,q_2})}_{L^2(\p\M)} \\ \leq
C \sum_{k\ge 1}\frac{\lambda_{k,q_1}(\lambda_{k,q_1}+\lambda_{k,q_2})}{(\lambda_{k,q_2}-z)^{n+2}}\leq C\epsilon.
\end{multline}
Therefore, we find
\begin{equation}\label{II.2.15}
\norm{\mathcal{I}_2(z)h}_{H^{1/2}(\p\M)}\leq C\norm{h}_{L^2(\p\M)}\epsilon.
\end{equation}
For $\mathcal{I}_3(z)h$, we have
\begin{align*}
\norm{\mathcal{I}_3(z)h}_{H^{1/2}(\p\M)} &\leq C\norm{h}_{L^2(\p\M)}\sum_{k\geq 1}
\frac{1}{\lambda_{k,q_2}^{n+1}} \norm{\partial _\nu \phi_{k,q_1}-\partial _\nu \phi_{k,q_2}}_{H^{1/2}(\p\M)}
\\
&\leq C\norm{h}_{L^2(\p\M)}\sum_{k\geq 1} \frac{1}{ k^{
\frac{2(n+1)}{n} } } \norm{\partial _\nu \phi _{k,q_2}-\partial
_\nu \phi_{k,q_1}}_{H^{1/2}(\p\M)}
\\
&\leq C\norm{h}_{L^2(\p\M)}\sum_{k\geq 1} \frac{1}{ k^{\frac{2r}{n} } } \norm{\partial _\nu \phi _{k,q_2}-\partial
_\nu \phi _{k,q_1}}_{H^{1/2}(\p\M)}.
\end{align*}
Therefore
\begin{equation}\label{II.2.16}
\norm{\mathcal{I}_3(z)h}_{H^{1/2}(\p\M)} \leq C\epsilon\norm{h}_{L^2(\p\M)}.
\end{equation}
The conclusion follows then from a combination of (\ref{II.2.16}),
(\ref{II.2.15}) and (\ref{II.2.11}).
\end{Demo}
\begin{Demo}{ of Theorem \ref{Th3}} From (\ref{II.2.9}) and Lemma \ref{LII.2.2}, we obtain
$$
\norm{ P^{(j)}(0)} _s \leq C\para{\abs{z}^{-j-\frac{1-2s}{4}}+\abs{z}^{n-j+1}\epsilon^{\mu_1}},
$$
and then
$$
\norm{ P^{(j)}(0)} _s \leq C \para{\abs{z}^{-\frac{1-2s}{4}}+\abs{z}^{n+1}\epsilon^{\mu_1}}, \quad \textrm{ if } \abs{z}\geq 1.
$$
In particular
\begin{equation}\label{II.2.17}
\norm{ P^{(j)}(0)} _s \leq C \min_{\rho \geq
1}\para{\rho^{-\frac{1-2s}{4}}+\rho ^{n+1}\epsilon^{\mu_1}}=C\epsilon^{\mu_2}
\end{equation}
where $\mu_2\in(0,1)$. Let $\mathscr{R}_{\g,q}$ be defined as in Lemma \ref{LII.2.3}. We
can proceed as in the proof of Lemma \ref{LII.2.4} to prove
\begin{equation}\label{II.2.18}
 \norm{\mathscr{R}_{\g,q_1}-\mathscr{R}_{\g,q_2}}_{\mathcal{L}\para{\mathcal{H}_1,\mathcal{H}_2}}\leq C\epsilon^{\mu_3} .
\end{equation}
From identity (\ref{II.2.1}), estimates (\ref{II.2.18}) and
(\ref{II.2.17}), we deduce
$$
\norm{\Lambda^\sharp_{\g,q_1}-\Lambda^\sharp_{\g,q_2}}_{\mathcal{L}\para{\mathcal{H}_1,\mathcal{H}_2}}\leq C\epsilon^{\mu_4},
$$
provided that $\epsilon$ is sufficiently small. To finish, we only need to remark that the traces of the geometrical optics solutions constructed in section
\ref{sec:GOSol} in fact satisfy
     $$ u_{1}|_{\p \M},u_{2}|_{\p \M} \in \mathcal{H}_{1} $$
so that in the proof of Theorem \ref{Th1} the right-hand side of \eqref{1.13} may be replaced by
     $$ \norm{\Lambda^\sharp_{\g,q_1}-\Lambda^\sharp_{\g,q_2}}_{\mathcal{L}\para{\mathcal{H}_1,\mathcal{H}_2}}. $$
This completes the proof of Theorem \ref{Th3}.
\end{Demo}

\subsection*{Acknowledgements} Part of this work was done while one of the author was visiting Tunisia. David DSF wishes to
express his gratitude for the hospitality of the Facult\' e des Sciences de Bizerte.


\begin{thebibliography}{99}
%
\bibitem{[AKKLT]}{M.~Anderson, A.~Katsuda,Y.~Kurylev,M.~Lassas,M.~Taylor, }{\it Boundary regularity for the Ricci equation, geometric convergence
and Gel'fand's inverse boundary problem}, Inventiones Math. 158 (2004), 261-321.

\bibitem{[AS]}{G.~Alessandrini and J.~Sylvester, }{\it Stability for multidimensional inverse spectral problem}, Commun. PDE, 15, (5) (1990), 711-736.
%
\bibitem{[AB]}{S.~A.~Avdonin, M.~I.~Belishev, }{\it Dynamical inverse problem for the Schr\"odinger equation (BC-method)}.
Uraltseva, N.N.(ed.), Proceedings of the St. Petersburg Mathematical Society. Vol. X. Transl. from the Russian by Tamara Rozhkovskaya. Providence,
RI: American Mathematical Society (AMS). Translations. Series 2. American Mathematical Society 214, 1-14 (2005); translation from Tr. St-Peterbg. Mat.
Obshch. 10, 3-17 (2004).
%
\bibitem{[1]}{M.~Belishev, }{\it Boundary control in reconstruction of manifolds and metrics (BC
method)}, Inverse Problems 13 (1997), R1 - R45.
%
\bibitem{[BK]}{M.~Belishev, Y.V.~Kurylev, }{\it To the reconstruction of a Riemannian manifold via its spectral data (BC- method)}.
Commun. Partial Differ. Equations 17, No.5-6, 767-804 (1992).
%
\bibitem{[Bellassoued]}{M.~Bellassoued, }{\it Uniqueness and stability in determining the speed of
propagation of second-order hyperbolic equation with variable
coefficients}, Applicable Analysis 83 (2004), 983-1014.
%
%
\bibitem{[BCY]}{M.~Bellassoued, M.~Choulli, M.~Yamamoto, }{\it Stability estimate for an inverse wave equation and a multidimensional Borg-Levinson theorem}, J. Diff. Equat. 247 (2) (2009) 465-494.
%
\bibitem{[Bel-Choul]}{M.~Bellassoued, M.~Choulli, }{\it Stability estimate for an inverse problem for the magnetic Schr\"odinger equation from the
Dirichlet-to-Neumann map}, J. Funct. Anal. 258, No. 1, 161-195 (2010).
%
\bibitem{[Bell-Jel-Yama1]}{M.~Bellassoued, D.~Jellali, M.~Yamamoto, }{\it Lipschitz stability for a hyperbolic inverse problem by
finite local boundary data,} Applicable Analysis 85 (2006),
1219-1243.
%
\bibitem{[Bell-Jel-Yama2]}{M.~Bellassoued, D.~Jellali, M.~Yamamoto, }{\it Stability Estimate for the hyperbolic inverse boundary value problem
by local Dirichlet-to-Neumann map}, J. Math. Anal. Appl. 343, No. 2, 1036-1046, (2008).
%
%
\bibitem{[Bukhgeim-Uhlmann]}{A.~L.~Bukhgeim and G.~Uhlmann, }{\it Recovering a potential from Cauchy data}, Comm. Partial Diff.
Equations 27 (2) (2002), 653-668.
%
\bibitem{[Calderon]}{A.~P.~Calder\'on, }{\it On an inverse boundary value problem,} in Seminar on
Numerical Analysis and its Applications to Continuum Physics, Rio
de Janeiro, (1988), 65-73.
%
\bibitem{[CM]}{F.~Cardoso and R.~Mendoza, }{\it On the hyperbolic Dirichlet-to-Neumann functional,} Comm.
Partial Diff. Equations 21 (1996), 1235-1252.
%
%
\bibitem{[CN]}{J.~Cheng and G.~Nakamura, }{\it Stability for the inverse potential problem by finite
measurements on the boundary,} Inverse Problems 17 (2001), 273-280.
%
\bibitem{[CY]}{J.~Cheng and M.~Yamamoto, }{\it The global uniqueness for determining two convection coefficients from Dirichlet-to-Neumann map in
two dimensions,} Inverse Problems 16, (3) (2000), L31-L38.
%
\bibitem{[Ch]}{M.~Choulli, }{\it Une Introduction aux Probl\`emes Inverses Elliptiques et Paraboliques},
Math\'ematiques et Applications, Vol. 65, Springer-Verlag, Berlin, 2009.
%
\bibitem{[DKSU]}{D.~Dos Santos Ferreira, C.~E.~Kenig, M.~Salo, and G.~Uhlmann, }{\it Limiting Carleman weihts and anisotropic inverse problems},
Inventiones Math. 178 (2009), 119-171.
%
\bibitem{[Eskin1]}{G.~Eskin, }{\it A new approach to hyperbolic inverse problems}, arXiv:math/0505452v3 [math.AP], 2006.
%
\bibitem{[Eskin2]}{G.~Eskin, }{\it Inverse hyperbolic problems with time-dependent coefficients}, arXiv:math/0508161v2 [math.AP], 2006.
%
\bibitem{[Eskin3]}{G.~Eskin, }{\it Global uniqueness in the inverse scattering problem for the Schr\"odinger operator with external
Yang-Mills potentials}, Comm. Math. Phys. 222 (2001), no. 3,
503-531.
%
\bibitem{[Eskin4]}{G.~Eskin, }{\it Inverse scattering problem in anisotropic media}, Comm. Math. Phys. 199 (1998), no. 2, 471-491.
%
%
\bibitem{[Hebey]}{E.~Hebey, }{\it Sobolev spaces on Riemannian manifolds}. Lecture Notes in Mathematics. 1635. Berlin: Springer. (1996).
\bibitem{[Hech-Wang]}{H.Hech-J-N.Wang, }{\it Stablity estimates for the inverse boundary value problem by partial Cauchy data},
 Inverse Problems, 22 (2006), 1787-1796.
%
\bibitem{[Ikawa]}{M.~Ikawa, }{\it Hyperbolic partial differential equations and wave phenomena}, Translations of Mathematical Monographs.
189. Providence, RI: American Mathematical Society (AMS). xxi, 190 p. (2000).
%
\bibitem{[Isakov1]}{V.~Isakov, }{\it An inverse hyperbolic problem with many boundary
measurements}, Comm. Part. Dif. Equations 16 (1991), 1183-1195.
%
\bibitem{[I1]}{V.~Isakov, }{\it Inverse Problems for Partial Differential Equations},
Springer-Verlag, Berlin, (1998).
%
\bibitem{[IS]}{V.~Isakov and Z.~Sun, }{\it Stability estimates for hyperbolic inverse problems with
local boundary data}, Inverse Problems 8 (1992), 193-206.
%
\bibitem{[Jost]}{J.~Jost, }{\it Riemannian Geometry and Geometric Analysis},  Universitext, Springer,
New York, 1995,
%
\bibitem{[KKL]}{A.~Katchalov, Y.~Kurylev and M.~Lassas, }{\it Inverse Boundary Spectral Problems,} Chapman \& Hall/CRC,
Boca Raton, (2001).
%
\bibitem{[KL]}{Y.V.~Kurylev and M.~Lassas, }{\it Hyperbolic inverse problem with data on a part of the
boundary,} in "Differential Equations and Mathematical Physics",
AMS/IP Stud. Adv. Math. 16, Amer. Math. Soc., Providence, (2000),
259-272.
%
%
\bibitem{[Lions-Magenes]}{J.-L.~Lions and E.~Magenes, }{\it Non-homogenous Boundary Value Problems and Applications},
Volumes I and II, Springer-Verlag, Berlin, (1972).
%
\bibitem{[PU]}{L.~Pestov and G.~Uhlmann, }{\it Two dimensional compact simple manifolds are boundary distance rigid},
Annals of Math., 161(2005), 1089-1106.
%
\bibitem{[Rachele]}{L.~Rachele, }{\it Uniqueness in inverse problems for elastic media with residual stress,} Comm. Partial Diff.
Equations 28 (2003), 1787-1806.
%
\bibitem{[Rakesh1]}{Rakesh, }{\it Reconstruction for an inverse problem for the wave equation
with constant velocity}, Inverse Problems 6 (1990), 91-98.
%
\bibitem{[Rakesh-Symes]}{Rakesh and W.~Symes, }{\it Uniqueness for an inverse problems for the wave equation},
Comm. Partial Diff. Equations 13 (1988), 87-96.
%
\bibitem{[Ramm-Sjostrand]}{A.~Ramm and J.~Sj\"ostrand, }{\it An inverse problem of the wave equation}, Math. Z. 206
(1991), 119-130.
%
\bibitem{[Salo]}{M.~Salo, }{\it Semiclassical Pseudodifferential Calculus and the Reconstruction
of a Magnetic Field}, arXiv:math/0602290v1 [math.AP] 14 Feb 2006.
%
\bibitem{[Sh]}{V.~Sharafutdinov, }{\it Integral Geometry of Tensor Fields}. VSP, Utrecht, the Netherlands, 1994.
%
\bibitem{[SU]}{P.~Stefanov and G.~Uhlmann, }{\it Stability estimates for the hyperbolic Dirichlet-to-Neumann
map in anisotropic media,} J. Functional Anal. 154 (1998),
330-358.
%
\bibitem{[SU1]}{P.~Stefanov and G.~Uhlmann, }{\it Stability estimates for the X-ray transform of tensor fields and boundary rigidity,}
Duke Math. J. 123(2004), 445-467.
%
\bibitem{[SU2]}{P.~Stefanov and G.~Uhlmann, }{\it Stable Determination of Generic Simple Metrics from the Hyperbolic Dirichlet-to-Neumann Map},
International Math. Research Notices, 17 (2005), 1047-1061.
%
\bibitem{[Sun]}{Z.~Sun, }{\it On continous dependence for an inverse initial boundary value
problem for the wave equation}, J. Math. Anal. App. 150 (1990),
188-204.
%
%
\bibitem{[Uhlmann]}{G.~Uhlmann, }{\it Inverse boundary value problems and applications}, Ast\'erisque 207, (1992), 153-221.

\end{thebibliography}
\end{document}